\documentclass{aims} % Use the aims.cls file to compile your paper
\usepackage{amsmath}
\usepackage{paralist}
\usepackage[misc]{ifsym}
\usepackage{epsfig} %For pictures: screened artwork should be set up with an 85 or 100 line screen
\usepackage{epstopdf} %This is to transfer .eps figure to .pdf figure; please compile your article using PDFLaTex or PDFTeXify.
\usepackage[colorlinks=true]{hyperref}
\hypersetup{urlcolor=blue, citecolor=red}
\allowdisplaybreaks

\def\currentvolume{X}
 \def\currentissue{X}
  \def\currentyear{200X}
   \def\currentmonth{XX}
    \def\ppages{X-XX}
     \def\DOI{10.3934/xx.xxxxxxx}

\newtheorem{theorem}{Theorem}[section]
\newtheorem{corollary}[theorem]{Corollary}

\newtheorem{lemma}[theorem]{Lemma}
\newtheorem{proposition}[theorem]{Proposition}

\theoremstyle{definition}

\newtheorem{assumption}[theorem]{Assumption}

\newcommand{\ep}{\varepsilon}
\newcommand{\eps}[1]{{#1}_{\varepsilon}}

\usepackage{amsfonts,amssymb,amsthm}
\usepackage{graphicx}
\usepackage{subcaption}
\newcommand{\rd}{\mathrm{d}}
\newcommand{\Kn}{\mathsf{Kn}}
\newcommand{\St}{\mathsf{St}}
\DeclareMathOperator*{\argmin}{arg\,min}
\newcommand{\Tr}{\text{t}}
\newcommand{\Sub}{\text{s}}
\newcommand{\tr}{\text{t}}
\newcommand{\sub}{\text{s}}

\title[Reconstruction of reflection coefficient]
{Stability of the reconstruction of the heat reflection coefficient in the phonon transport equation} % Only the first word and proper nouns should be capitalized.

\author[Peiyi Chen, Irene M. Gamba, Qin Li and Anjali Nair]{}


\thanks{$^*$Corresponding author: Peiyi Chen}

\begin{document}
\maketitle

% Enter the first author's name and email address; email addresses are required for each author.
% Use footnote notations to indicate address and affiliations with commas between numbers if more than one address applies; see below for examples.
\centerline{\scshape
Peiyi Chen$^{{\href{pchen345@wisc.edu}{\textrm{\Letter}}}*1}$,
Irene M. Gamba$^{{\href{gamba@math.utexas.edu}{\textrm{\Letter}}}2}$,
Qin Li$^{{\href{qinli@math.wisc.edu}{\textrm{\Letter}}}1}$,
and Anjali Nair$^{{\href{anjalinair@uchicago.edu}{\textrm{\Letter}}}3}$}

\medskip

{\footnotesize
% Enter the full affiliation and country name:
% Do not insert commas or periods at the end of lines.
 \centerline{$^1$Department of Mathematics, University of Wisconsin--Madison, USA}
} % Do not forget to end {\footnotesize with the sign }

\medskip

{\footnotesize
 % Enter the full affiliation and country name:
 \centerline{$^2$Department of Mathematics and Oden Institute for Computational Engineering and Sciences,
The University of Texas at Austin, USA}
}

\medskip

{\footnotesize
 % Enter the full affiliation and country name:
 \centerline{$^3$Department of Mathematical Sciences, University of Delaware, USA}
}

\bigskip

% The name of the handling editor will be entered by AIMS production staff.
% "Communicated by Handling Editor" is not needed for special issue.
 \centerline{(Communicated by Handling Editor)}

%%%%%%%%%%%%%%%%%%%%%%%%%%%%%%%%%%%%%%%%%%%%%%%%%%%%%%%
%             5. ABSTRACT
%%%%%%%%%%%%%%%%%%%%%%%%%%%%%%%%%%%%%%%%%%%%%%%%%%%%%%%

\begin{abstract}
The reflection coefficient is an important thermal property of materials, especially at the nanoscale, and determining this property requires solving an inverse problem based on macroscopic temperature measurements.
In this manuscript, we investigate the stability of this inverse problem to infer the reflection coefficient in the phonon transport equation. We show that the problem becomes ill-posed as the system transitions from the ballistic to the diffusive regime, characterized by the Knudsen number converging to zero. Such a stability estimate clarifies the discrepancy observed in previous studies on the well-posedness of this inverse problem. Furthermore, we quantify the rate at which the stability deteriorates with respect to the Knudsen number and confirm the theoretical result with numerical evidence.
\end{abstract}

%%%%%%%%%%%%%%%%%%%%%%%%%%%%%%%%%%%%%%%%%%%%%%%%%%%%%%
%                   6. BODY
%%%%%%%%%%%%%%%%%%%%%%%%%%%%%%%%%%%%%%%%%%%%%%%%%%%%%%

% Only the first word and proper nouns of section titles should be capitalized.
% The title of section 1:
\section{Introduction}

Material design for desired thermal properties~\cite{Wang2011-nano_silicon}, mathematically, is related to solving a PDE-based inverse problem. Usually, it comes down to using the temperature response of the material to infer its thermal property, represented by certain parameters in a PDE that characterizes heat propagation.

To describe heat propagation, the classical heat equation requires the validity of Fourier's law, which is no longer valid in the microscale, or nanoscale~\cite{Regner_FDTR_heat_conduct_2013, Johnson_2013_nondiffusive_measure}. Instead, thermal propagation is described by the following Phonon Transport Equation (PTE) derived from the first principle~\cite{GC_1996_nano_BTE, Majumdar_1993_Microscale, hua2017experimental}, %at this scale.
\begin{equation}\label{eq:nonlinear_PTE_F-eq}
\partial_t F + \mu\nu(\omega)\partial_x F = \frac{F^\ast - F}{\tau(\omega)}\,,
\end{equation}
where $F(t,x,\mu,\omega)$ denotes the energy distribution function of phonon at time $t$, position $x$ with direction cosine $\mu$ and phonon frequency $\omega$; $\nu(\omega)$, $\tau(\omega)$ denote the group velocity and the phonon relaxation time, respectively, and they both depend on the phonon frequency. Then $F^\ast(t,x,\mu,\omega)=\hbar\omega D(\omega)\left(e^{\frac{\hbar\omega}{k_BT(t,x)}}-1\right)^{-1}$ is the local equilibrium given by the Bose-Einstein distribution. The more precise definition is in~\eqref{eq:energy_conserve_nl_PTE}. This is a reduced model assuming plane symmetry, hence $\partial_y F, \partial_z F$ are dropped.

Usually in a lab environment, material properties (such as $\nu$ and $\tau$) are typically hard to measure directly, and it is a standard practice to conduct experiments to measure macroscopic quantities and use them to infer these quantities, naturally introducing an inverse problem.

\begin{figure}[htbp]
    \centering
    \includegraphics[width=12cm]{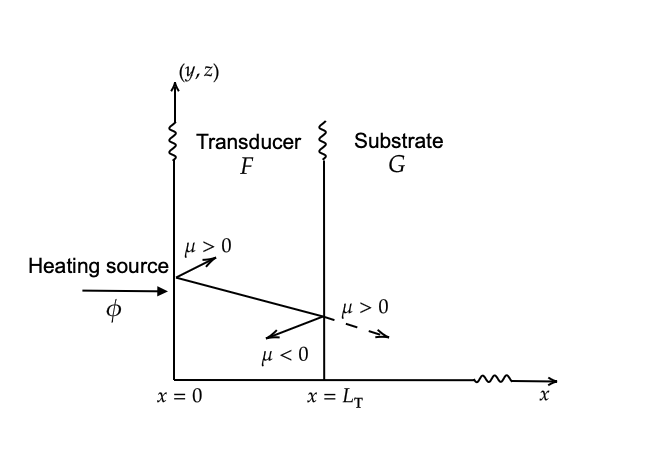}
  \caption{Schematics of the experimental setup~\cite{hua2017experimental}. Here the angle corresponds to the directional cosine $\mu$.}
  \label{fig:1}
\end{figure}

We study one such experiment in this paper that is deployed for analyzing reflection coefficients. More specifically, two types of materials (one transducer and one substrate) are placed adjacent to each other, 
and we are tasked to understand the heat conductance at the interface. To do so, a heat source is injected through the surface of the transducer, and one can measure its temperature response as a function of time. This data pair can potentially be used to infer the reflection coefficient at the interface, see Figure~\ref{fig:1}.

Mathematically, thermal transport 
is modeled as a system of two phonon transport equations coupled by the boundary conditions encoding the transmission and reflection coefficients. The phonon density within the transducer is denoted by $F$, solving
\begin{equation}\label{eq:nonlinear-PTE-F}
\left\{
\begin{aligned}
    \partial_t F + \mu\nu_{\Tr}(\omega)\partial_x F =& \frac{F^\ast - F}{\tau_{\Tr}(\omega)}\,, & 
    x\in [0,L_{\Tr}] \\
    F(t,x=0,\mu,\omega) =& F^\ast_0+ \phi(t,\mu,\omega) \,, & \mu\in[0,1]\\
    F(t,x=L_{\Tr},\mu,\omega) =& 
    F^\ast_0 + \eta_{\tr}(\omega)(F(t,x=L_{\Tr},-\mu,\omega)-F^\ast_0) \\
    & \quad\ + \zeta_{\tr}(\omega)(G(t,x=L_{\Tr},\mu,\omega)-F^\ast_0)\,, & \mu\in[-1,0]%\\
    %F(t=0,x,\mu,\omega) =& F^\ast_0(\omega)\,.
\end{aligned}
\right. 
\end{equation}
and the phonon density within the substrate, $G$, solves
\begin{equation}\label{eq:nonlinear-PTE-G}
\left\{
\begin{aligned}
    \partial_t G + \mu\nu_{\Sub}(\omega)\partial_x G = & \frac{F^\ast - G}{\tau_{\Sub}(\omega)}\,, & x\in [L_{\Tr},L_{\Tr}+L_{\Sub}] \\
    G(t,x=L_{\Tr},\mu,\omega) = & F^\ast_0 + \eta_{\sub}(\omega)(G(t,x=L_{\Tr},-\mu,\omega)-F^\ast_0) \\
    & \quad\ + \zeta_{\sub}(\omega)(F(t,x=L_{\Tr},\mu,\omega)-F^\ast_0) \,, & \mu\in[0,1]\\
    G(t,x=L_{\Tr}+L_{\Sub},\mu,\omega) = & F^\ast_0 \,, & \mu\in[-1,0]%\\
    %G(t=0,x,\mu,\omega) = & F^\ast_0\,.
\end{aligned}
\right. 
\end{equation}
where the subscripts $\Tr, \Sub$ stand for the transducer and the substrate, respectively.
The coefficients $\eta_{\tr}(\omega), \eta_{\sub}(\omega)$ and $\zeta_{\tr}(\omega), \zeta_{\sub}(\omega)$ are the reflection and transmission coefficients that are to be inferred from the measurements of the temperature response. 
All these coefficients take the value between $0$ and $1$. By energy conservation of the system across the interface~\cite{hua2017experimental}, and assuming transmission and reflection processes independent of phonon frequency, we obtain 
\begin{equation}\label{eq:tr_sub_couple}
\begin{aligned}
1 = \eta_{\tr}(\omega) + \frac{\nu_{\Sub}(\omega)}{\nu_{\Tr}(\omega)} \zeta_{\sub}(\omega)\,,\quad
% \nu_{\Tr} (\omega) &= \nu_{\Tr}(\omega)\eta_{\tr}(\omega) + \nu_{\Sub}(\omega)\zeta_{\sub}(\omega)\,,\\
1 = \frac{\nu_{\Tr}(\omega)}{\nu_{\Sub}(\omega)}\zeta_{\tr}(\omega) + \eta_{\sub}(\omega)\,.
% \nu_{\Sub}(\omega) &= \nu_{\Tr}(\omega)\zeta_{\tr}(\omega) +  \nu_{\Sub}(\omega) \eta_{\sub}(\omega)\,.
\end{aligned}
\end{equation}
In addition, the following compatibility condition,
\begin{equation}\label{eq:eta_1_compatible_condition}
\eta_{\tr}(\omega) + c \zeta_{\tr}(\omega) =1\,, 
\ \text{for some positive constant } c=O(1)\,, 
\ \textit{(detailed balance condition)}
\end{equation}
further reduces the system with the unknown pair $(\eta_{\tr},\zeta_{\tr})$ into a scalar initial boundary value problem for only $\eta_{\tr}(\omega)$.  This requirement arises from physics~\cite{hua2017experimental}, which ensures the existence of the nontrivial steady state equilibrium for the $(F, G)$ vector-valued system, and more details will be given for the linearized model in section~\ref{subsec:full_model}. 
The incoming boundary data 
is of Dirichlet type, prescribed by $F(x=0) = F_0^\ast + \phi$, for $\phi$ externally given. The initial data for the coupled system are prepared as the global equilibrium at room temperature, given by $(F,G)|_{t=0}=(F_0^\ast, F_0^\ast)$.

The experiments were first studied in~\cite{hua2017experimental} with the Time-Domain ThermoReflectance (TDTR) mechanism. Follow-up experiments and data fitting were conducted in~\cite{forghani2018reconstruction, Hua_Lindsay_general_Fourier_2019}. Mathematically, the problem was modeled and investigated in~\cite{gamba2022reconstructing} where the authors formulated the reconstruction of the reflection coefficient at the boundary as a PDE-constrained optimization and numerically conducted the inverse problem by stochastic gradient descent (SGD). The first attempt at theoretical inversion was completed in~\cite{sun2022unique}, where the authors applied the singular decomposition technique and gave an explicit reconstruction formula.

Despite all these efforts, mathematical investigation and the physics literature lead to inconsistent conclusions. Indeed, in~\cite{sun2022unique} the authors claimed well-posedness for the reconstruction, but~\cite{Minnich_MFP_2012} reports that the simple data-fitting has to be modulated with relaxation terms for a meaningful reconstruction. Though it is not uncommon to have misaligned mathematics and physics results, the drastic disparity in this formulation is alarming. The current paper %forms one attempt to 
addresses and understands this disparity by the mathematical framework described as follows.

We find that a big component of the inconsistency in previous studies is due to the fact the incorrect regime the coupled system~\eqref{eq:nonlinear-PTE-F}--\eqref{eq:nonlinear-PTE-G} is in for analysis. Indeed, in~\cite{gamba2022reconstructing, sun2022unique}, the authors started off assuming the non-dimensional system with hyperbolic structure and directly dived in with the reconstruction process. As a typical hyperbolic system, information propagates along characteristics, and it is expected that one can trace back the information for the inversion and deduce a well-posed inverse problem. One phenomenon, however, that gets neglected in those papers is the long-time asymptotic heat equation limit. Indeed, in the long-time large space regime, the heat equation is regarded as an accurate model for characterizing heat propagation. As a typical parabolic equation, the heat equation is well-known to behave poorly for inverse problems. This means that in this limit, the well-posedness for the PTE inverse problem, as proved in~\cite{sun2022unique}, should have its nice features deteriorated.

This is exactly our plan. To mathematically describe this deterioration phenomenon, we deploy $\ep$ to denote the Knudsen number that codes the scaling of the system. Through investigating parameters used in lab experiments, and performing the non-dimensionalization for the system, we analyze the associated regimes, and trace the well-posedness Lipschitz constant using $\ep$. We will discover that some experimental configuration actually lands us in the macroscopic regime where the heat equation is prevalent. When parabolic features are dominant, ill-conditioning is unavoidable. Mathematically, denoting $\mathcal{M}^\ep(\eta_{\tr})$ the measurement induced by the reflection coefficient $\eta_{\tr}$, we are to prove:
\begin{equation}
\|\mathcal{M}^\ep(\eta_1)-\mathcal{M}^\ep(\eta_2)\|\lesssim \|\eta_1-\eta_2\|e^{-\tilde{c}/\ep}\,,\quad \text{for some } \tilde{c}>0\,,
\end{equation}
for different transducer reflection coefficients $\eta_1, \eta_2$ in suitable norms. That is, in the diffusion limit ($\eps\to0^+$), the Lipschitz constant converges to zero, and the temperature measurement becomes indistinguishable, making the inversion ill-posed. A more precise statement will be presented in Theorem~\ref{thm:stab_estim}.

It is not new research to analyze the impact of PDEs' scaling on inverse problems, especially in the kinetic framework~\cite{JNWang_stable_transport_1999, Gunther_EIT_unstable_2009, JNWang_2013_increase_stable, li2015diffusion, Zhao_Zhong_2018_RTE_Instability,  Lai_Li_Uhlmann_RTE_2019, cheng2011recovering, 
Albi_Herty_Jorres_Pareschi_2014_AP_control, Herty_Visconti_2019_inverse_EnKF, Albi_Herty_Pareschi_2015_opinion, Herty_Puppo_Roncoroni_Visconti_2020_traffic, Doumic_Moireau_2026}. The main technique is the singular decomposition, a method widely used for inverse problems in kinetic equation setting~\cite{Choulli_Stefanov_reconstruct_1996, Choulli_Stefanov_singular, SR_Arridge_opto_review_1999, Bal_transport_review_2009, Bal_Jollivet_stable_transport_2010, Bal_average_2011, chen_Li_Wang_RTE_2018, Li_Sun_singular_2020}.

The paper is organized as follows. In Section~\ref{sec:PTE_model}, we linearize the PTE and introduce the diffusive scaling, under which the Strouhal number and the Knudsen number are comparable and assumed to be small. The complete PTE system is summarized in~\eqref{eq:forward_PTE_f_eps_couple} and~\eqref{eq:forward_PTE_g_eps_couple}. The inverse problem is then formulated in Section~\ref{sec:inverse_results}, followed by the statement of our main result in Section~\ref{subsec:stability_main_thm}. Section~\ref{sec:pf_stability_thm} is devoted to the proof of Theorem~\ref{them:measurem_decomp}, where we employ the singular decomposition technique and conduct the estimates. Finally, numerical results are presented in Section~\ref{sec:numerics}.

% The title of section 2:
\section{Phonon Transport Equation and non-dimensionalization}\label{sec:PTE_model}

We summarize in section~\ref{subsec:setup} the main setup and parameters in lab experiments as reported in~\cite{hua2014transport, hua2017experimental, forghani2018reconstruction, Peraud_Hadjiconstantinou_2016_PTE}. The measurement suggests that the temperature variation is small enough, making linearization around room temperature a valid assumption. To match the scaling in experiments, we also conduct non-dimensionalization. These procedures are discussed in section~\ref{subsec:nondimensionalize}.

\subsection{Scalings from experimental setups}\label{subsec:setup}
We review the experimental setups reported in~\cite{hua2014transport, hua2017experimental, forghani2018reconstruction, Peraud_Hadjiconstantinou_2016_PTE}. Two parallel materials are placed adjacent to each other, with the thin metal (aluminum) absorbing heat, and attached to a much thicker substrate (silicon). A laser source is shined on the metal. The laser frequency can be adjusted to trigger the response of phonons at different frequencies. Sensors that measure temperature response are placed at the surface of the metal to detect temperature increases as a function of time.

The parameters used in the experiments are summarized in the following table~\ref{tab:parameter}.
\begin{table}[htbp]
    \centering
    \begin{tabular}{l c c}
    \textit{Parameters} & Notation & Units \\
    \hline
    Thickness of the transducer (Aluminum) & $L_{\Tr}$
    %$\mathtt{L}$
    & $10^{-7}-10^{-5} \, \mathrm{m}$\\
    Thickness of the substrate (Silicon) & & - \\
    Sensor measuring time %gaps
    & $\mathtt{T}$ & $10^{-8}-10^{-5}\,\mathrm{s}$\\
    Phonon frequency & $\omega$ & $0-100\,\mathrm{THz}$ \\
    Room temperature & $T_0$ & $300 \,\mathrm{K}$\\
    Range of temperature increase & $\Delta T$ & $0-1\, \mathrm{K}$\\
    Characteristic bulk velocity & $L_{\Tr}/\mathtt{T}$ & $10\, \mathrm{m}\cdot\mathrm{s}^{-1}$\\
    Characteristic group velocity & $\bar{\nu}$ & $10^{3} \,\mathrm{m}\cdot\mathrm{s}^{-1}$ \\
    Group velocity of the transducer (Aluminum) & $\nu_{\Tr}$ & $5000\,\mathrm{m}\cdot\mathrm{s}^{-1}$\\
    Group velocity of the substrate (Silicon) & $\nu_{\Sub}$ & $2600\,\mathrm{m}\cdot\mathrm{s}^{-1}$\\
    Characteristic relaxation time & $\bar{\tau}$ & $10 \,\mathrm{ps}$ \\
    Relaxation time of the transducer (Aluminum) & $\tau_{\Tr}$ & $10\,\mathrm{ps}$\\
    Relaxation time of the substrate (Silicon) & $\tau_{\Sub}$ & $40\,\mathrm{ps}$\\
    Characteristic mean free path (MFP) & $\bar{\Lambda}$ & $10 \,\mathrm{nm}$\\
    \hline
    \textit{Physical constants} \\
    \hline
    Reduced Planck constant & $\hbar$ & $10^{-34} \mathrm{m}^{2} \cdot \mathrm{kg} \cdot \mathrm{s}^{-1}$ \\
    Boltzmann constant & $k_B$ & $10^{-23}\mathrm{J}\cdot \mathrm{K}^{-1}$
    \end{tabular}
    \caption{Parameters and physical constants}
    \label{tab:parameter}
\end{table}

As is clear in the data above, the temperature variation ($0-1\, \mathrm{K}$) is much smaller compared to room temperature, so the system can be safely considered in the linearized regime. We therefore conduct linearization in subsection~\ref{subsec:linearize}. Also, multiple quantities that share units span a wide range of scaling, demonstrating a multiscale feature of the problem. We conduct the rescaling and non-dimensionalization in subsection~\ref{subsec:nondimensionalize}. Since the substrate is significantly thicker than that of the transducer in the real experiment~\cite{hua2017experimental}, it takes much longer for its right-hand boundary condition to respond to the heat source and is regarded as of infinite length. 
We summarize the final model in section~\ref{subsec:full_model}.

\subsection{Linearization}\label{subsec:linearize}

We are to perform a linearization of the system around the room temperature, denoted by $T_0$, by linearly expanding $F$, the solution to~\eqref{eq:nonlinear_PTE_F-eq}, using:
\begin{equation}
f(t,x,\mu,\omega):=F-F_0^\ast = F-F^\ast(T_0)\,, 
\end{equation}
also known as the deviational energy distribution function.
Since both the time differentiation and the transport term are linear, the linearization is straightforward. We examine the linearization of the collision term on the right.

Note that the equilibrium distribution is given by the Bose-Einstein distribution, with its temperature determined through energy conservation. This leads to:
\begin{equation}\label{eq:energy_conserve_nl_PTE}
F^\ast=F^\ast(t,x,\mu,\omega)=\frac{\hbar\omega D(\omega)}{e^{\frac{\hbar\omega}{k_B T(t,x)}}-1}\,,\quad\text{with}\quad \int\limits_{-1}^1\int\limits_{0}^\infty\frac{F-F^\ast}{\tau} \rd\omega\rd\mu=0\,.
\end{equation}
Here, $D(\omega)>0$ is the phonon density of states, $T$ is the temperature of the material that depends on time and space, $\hbar$ is the reduced Planck constant, and $k_B$ is the Boltzmann constant. Linearizing this operator means writing $
\mathcal{L}f=\partial_T F^\ast|_{T_0}(T(t,x)-T_0)=: \partial_T F^\ast|_{T_0} \Delta T(t,x)$. Denoting $\langle f \rangle := \int\limits_{-1}^1\int\limits_{0}^\infty f\, \rd\omega\rd\mu$, and enforcing energy conservation by setting $\langle \frac{1}{\tau(\omega)}(\mathcal{L}f - f) \rangle =0$, we have:
\begin{equation}\label{eq:DeltaT_def}
\begin{aligned}
  & \Delta T=  \Delta T(t,x)=\frac{\langle  f/\tau\rangle(t,x)}{C_\tau}\,,\\
  & \text{with}\quad C_\omega :=\partial_T F^\ast|_{T_0} = \frac{\left(F^\ast\mid_{T_0}\right)^2 e^{\hbar\omega/k_B T_0}}{D(\omega) k_B T_0^2}\,,
  \ C_\tau :=\langle C_\omega/\tau\rangle\,.
\end{aligned}
\end{equation}
This gives:
\begin{equation}\label{eq:linear_Lf}
\mathcal{L}f=\frac{C_\omega}{C_\tau}\langle f/\tau\rangle\,.
\end{equation}
As a summary, the linearized equation of $f$ reads:
\begin{equation}\label{eq:linear_PTE_f}
\partial_t f + \mu\nu(\omega)\partial_x f = \frac{1}{\tau(\omega)}(\mathcal{L}f - f)\,.
\end{equation}

\subsection{Non-dimensionalization} \label{subsec:nondimensionalize}
% Assume that the laser heat source has a characteristic time scale of $\mathtt{T}$, i.e.,  
% \begin{equation}
% f(t,x=0,\mu>0,\omega) = \phi^{\mathtt{T}}(t,\mu,\omega) = \phi\Big(\frac{t}{\mathtt{T}},\mu,\omega\Big)\,.
% \end{equation}
Consider the characteristic time scale $\mathtt{T}$, the characteristic spatial scale $\mathtt{L}$, the characteristic group velocity $\bar{\nu}$ and the characteristic relaxation time $\bar{\tau}$, we introduce the rescaled coordinates $\left(\Tilde{t}:=\dfrac{t}{\mathtt{T}},\Tilde{x}:=\dfrac{x}{\mathtt{L}}\right)$, the rescaled group velocity $\tilde{\nu}(\omega):=\dfrac{\nu(\omega)}{\bar{\nu}}$, and the rescaled relaxation time $\tilde{\tau}(\omega):=\dfrac{\tau(\omega)}{\bar{\tau}}$.
In these rescaled variables, $\tilde{f}:=f(\tilde t,\tilde x,\mu,\omega)$ solves:
\begin{equation}\label{eq:rescale_PTE_tildef}
\dfrac{\mathtt{L}/\mathtt{T}}{\bar{\nu}}\partial_{\tilde t} \tilde{f} +\mu \tilde{\nu}(\omega) \partial_{\tilde x} \tilde{f} = \frac{\mathtt{L}}{\bar{\tau}\bar{\nu}} \frac{1}{\tilde{\tau}(\omega)} (\mathcal{L}\tilde f - \tilde f)\,,
\end{equation}
here we note that $\frac{\mathtt{L}/\mathtt{T}}{\bar{\nu}}=:\St$ is the Strouhal number that represents the ratio between the characteristic bulk velocity $\mathtt{L}/\mathtt{T}$ and the characteristic phonon group velocity $\bar{\nu}$, and $\frac{\bar{\nu} \bar{\tau}}{\mathtt{L}}=:\Kn$ is the Knudsen number that represents the ratio between the characteristic phonon mean free path, $\bar{\Lambda}:=\bar{\nu}\bar{\tau}$, and the domain size. That is, we can drop all the tildes and rewrite the equation as follows.
\begin{equation}\label{eq:rescale_PTE_f}
\St \partial_{t} f +\mu \nu(\omega) \partial_{x} f = \frac{1}{\Kn} \frac{1}{\tau(\omega)} (\mathcal{L} f - f)\,.
\end{equation}

According to Table~\ref{tab:parameter}, the characteristic mean free path is in the nanometer range, i.e., $\bar{\Lambda}= 10^{-8}\,\mathrm{m}$ and the characteristic relaxation time is in the range of picoseconds, i.e., $\bar{\tau}=10^{-11}\,\mathrm{s}$, then we have that $\bar{\nu}=10^3\, \mathrm{m}\cdot\mathrm{s}^{-1}$. In the diffusive regime, thickness of the transducer $\mathtt{L}\approx 10^{-5}\,\mathrm{m}$, and we have that
\begin{equation*}
    \Kn=\frac{\bar{\nu}\bar{\tau}}{\mathtt{L}}\approx 10^{-3}\,.
\end{equation*}
In addition, in the diffusive regime, the time scale $\mathtt{T}$ is in microseconds, $\mathtt{T}\approx10^{-5}\,\mathrm{s}$, 
\begin{equation*}
    \St= \frac{\mathtt{L}/\mathtt{T}}{\bar{\nu}}\approx 10^{-3}\,.
\end{equation*}
This places us perfectly in the regime consistent with a diffusion scaling, which is $\eps:=\St=\Kn$ for some small $0<\eps\ll 1$.
% \begin{equation}
%     \overline{\Kn}=\frac{\langle C_\omega \Kn(\omega)/\tau\rangle}{C_\tau},\quad \overline{\St}=\frac{\langle C_\omega \St(\omega)/\tau\rangle}{C_\tau}\,,
% \end{equation}
% then
% \begin{equation}
%      \overline{\Kn}= \overline{\St}\,.
% \end{equation}
The asymptotic analysis of the steady state phonon transport equation in the diffusive regime and the associated boundary layer analysis for boundary conditions was done in~\cite{nair2022second}.
Investigations of appropriate numerical algorithms in the vanishing Knudsen number limit have been conducted extensively~\cite{Jin_Pareschi_Toscani_1998_AP,Lafitte_Samaey_2012_AP_diff, Dimarco_Pareschi_Samaey_2018_MC_transport, JIN2022110895}.
%We briefly recall this in the Appendix~\ref{subsec:boundary_analysis} as applied to the time-dependent problem. %Here, we will assume a diffusive scaling of the form 
   % \begin{equation}
   %     \mathtt{L}=\frac{1}{\eps},\quad 
   %     \mathtt{T}=\frac{1}{\eps^2},\quad \eps\ll 1\,.
   % \end{equation}
% To explicitly demonstrate the $\eps$-dependence, we define $g^\eps(t,x,\mu,\omega)=f(t,x,\mu,\omega)$ that solves:
% We rewrite~\eqref{eq:rescale_PTE_f} as follows, and explicitly add $\eps$ in its superscript to emphasize its dependence.
% \begin{equation}
% \eps \partial_t f^\eps + \mu\nu(\omega)\partial_x f^\eps = \frac{1}{\eps} \frac{1}{\tau(\omega)}(\mathcal{L} f^\eps-f^\eps)\,.
% \end{equation}
% We drop the tilde for notational convenience. This gives a PDE for $g^\eps$ as
% \begin{equation}\label{eqn:PTE_diffusive}
% \begin{aligned}
% \eps^2\partial_tg^\eps+\eps\mu\nu(\omega)\partial_x g^\eps &=\frac{1}{\tau(\omega)}(\mathcal{L}g^\eps-g^\eps)\,, \\
% g^\eps(t,x=0,\mu,\omega)&=\phi(t,\mu,\omega)\,, & \mu>0\\
% g^\eps(t,x=1,\mu,\omega)&=\eta(\omega) g^\eps(t,x=1,-\mu,\omega)\,, & \mu<0\\
% g^\eps(t=0,x,\mu,\omega)&=0\,. 
% \end{aligned}
% \end{equation}

\subsection{The full model}\label{subsec:full_model}
We put all the pieces together. This is to linearize~\eqref{eq:nonlinear-PTE-F} for the transducer and~\eqref{eq:nonlinear-PTE-G} for the substrate, both around $F^\ast_0$, and conduct non-dimensionalization as demonstrated in the previous section. 
We identify $L_{\sub}$ as the dimensionless substrate thickness, and denote
$\nu_{\tr}, \tau_{\tr}$ as the dimensionless group velocity and relaxation time of the transducer, while the substrate quantities are defined analogously. By introducing the dimensionless number $\ep$, solutions of the coupled system are denoted by $f^\ep$ and $g^\ep$ respectively. 
We also identify $f^\ep= g^\ep=C_{\omega}$ as the unique homogeneous equilibrium solution to the coupled system, since they solve the same linear equation $(\mathcal{L}-I)[\cdot] = 0$.

% Together, the linearization of the PTE system~\eqref{eq:nonlinear-PTE-F} and~\eqref{eq:nonlinear-PTE-G} under the diffusive scaling reads as follows 
% \begin{equation}\label{eq:forward_PTE_g_eps}
%     \left\{\begin{aligned}
%     \eps \partial_t f^{\eps} +\mu\nu(\omega)\partial_x f^{\eps} &=\frac{1}{\eps} \frac{1}{\tau(\omega)}(\mathcal{L}f^{\eps} - f^{\eps})\,, \\
%     f^{\eps}(t,x=0,\mu,\omega)&=\phi(t,\mu,\omega), & \mu>0\\
%     f^{\eps}(t,x=1,\mu,\omega)&=\eta(\omega) f^{\eps}(t,x=1,-\mu,\omega), & \mu<0\\
%     f^{\eps}(t=0,x,\mu,\omega)&=0\,.
%     \end{aligned}\right.
% \end{equation}
% Let $f^\eps$ be the solution in the transducer and $g^\eps$ be the solution in the substrate as indicated in Figure~\ref{fig:1}. They are connected at the boundary $x=1$ by the transmission coefficient~\cite{sun2022unique}.
Altogether, this leads to the following.
\begin{equation}\label{eq:forward_PTE_f_eps_couple}
    \left\{\begin{aligned}
    \ep \partial_t f^{\ep} +\mu\nu_{\tr}(\omega)\partial_x f^{\ep} &=\frac{1}{\ep} \frac{1}{\tau_{\tr}(\omega)}(\mathcal{L}f^{\ep} - f^{\ep})\,, & x \in[0,1]\\
    f^{\ep}(t,x=0,\mu,\omega)&=\phi(t,\mu,\omega), & \mu>0\\
    f^{\ep}(t,x=1,\mu,\omega)&=\eta_{\tr}(\omega) f^{\ep}(t,x=1,-\mu,\omega) + \zeta_{\tr}(\omega) g^\ep(t,x=1,\mu,\omega), & \mu<0\\
    f^{\ep}(t=0,x,\mu,\omega)&=0\,.
    \end{aligned}\right.
\end{equation}
\begin{equation}\label{eq:forward_PTE_g_eps_couple}
    \left\{\begin{aligned}
    \ep \partial_t g^{\ep} +\mu\nu_{\sub}(\omega)\partial_x g^{\ep} &=\frac{1}{\ep} \frac{1}{\tau_{\sub}(\omega)}(\mathcal{L}g^{\ep} - g^{\ep})\,, \hspace{1.5in} x\in[1, 1+L_{\sub}] \\
    g^{\ep}(t,x=1,\mu,\omega)&= \eta_{\sub}(\omega) g^{\ep}(t,x=1,-\mu,\omega) + \zeta_{\sub}(\omega)f^\ep (t,x=1,\mu,\omega) \,, \hspace{.2in} \mu>0\\
    g^{\ep}(t,x=1+L_{\text{s}},\mu,\omega)&= 0\,, \hspace{3in} \mu<0\\
    g^{\ep}(t=0,x,\mu,\omega)&=0\,.
    \end{aligned}\right.
\end{equation}
The boundary condition at $x=0$ is of Dirichlet type in the incoming direction for $f^\ep$, and that at $x=1$ is of reflective type at the interface, with the reflection index $\eta_{\tr}(\omega)$ denoting the portion of phonons that get bounced back.
The explicit calculations of how the four coefficients at the boundary $(\eta_{\tr}, \zeta_{\tr}, \eta_{\sub}, \zeta_{\sub})$ are reduced to one unknown, $\eta_{\tr}$, is conducted in Appendix~\ref{apx:reduction}. 

According to the typical group velocities in Table~\ref{tab:parameter}, $\nu_{\sub}/\nu_{\tr}$ is approximately a constant of $O(1)$. 
In the special case where $\nu_{\tr}(\omega)= \nu_{\sub}(\omega)$, the transmission coefficients and the reflection coefficients are related by $\eta_{\tr}=1-\zeta_{\sub}$ and $\eta_{\sub}=1-\zeta_{\tr}$, which is an assumption commonly seen for coupled systems~\cite{Gamba_Degond_2002_couple}. In our current setup with aluminum and silicon being the transducer and substrate, respectively, it can be specified that $\nu_{\sub}/\nu_{\tr} = \frac{1}{2}$ and $\tau_{\sub}/\tau_{\tr} = 4$.

Moreover, the compatibility condition~\eqref{eq:eta_1_compatible_condition} requires that no heat can transmit across the boundary when both systems are at equilibrium~\cite{hua2017experimental}, and without loss of generality, it is set that $c\equiv1$.
Indeed, when $f^\ep = g^\ep =  C_\omega$, 
the boundary conditions coincide with the following convex combination satisfied  by the transmission and reflection coefficients~\eqref{eq:eta_1_compatible_condition},
\begin{equation}
\begin{aligned}
    \eta_{\tr}(\omega) + \zeta_{\tr}(\omega) = 1\,. %C_\omega = C_\omega\,, & (\text{boundary condition for } f^{\eps})\\
    % & \left[\eta_{\sub}(\omega) + \zeta_{\sub}(\omega) \right]C_\omega = \left[\left(1 - \frac{\nu_{\tr}}{\nu_{\sub}}\zeta_{\tr}(\omega)\right) + \frac{\nu_{\tr}}{\nu_{\sub}} \left(1- \eta_{\tr}(\omega)\right) \right]C_\omega = C_\omega\,, %& (\text{boundary condition for } g^{\eps})
    \end{aligned}
\end{equation}
With the aforementioned simplifications, the system of substrate $g^\ep$ now reads
\begin{equation}\label{eq:forward_PTE_g_eps_simple}
    \left\{\begin{aligned}
    \ep \partial_t g^{\ep} + \mu\frac{\nu_{\tr}(\omega)}{2} \partial_x g^{\ep} =&\frac{1}{\ep} \frac{1}{4\tau_{\tr}(\omega)}(\mathcal{L}g^{\ep} - g^{\ep})\,, & x\in[1, 1+L_{\sub}] \\
    g^{\ep}(t,x=1,\mu,\omega) =& (2\eta_{\tr}(\omega) -1) g^{\ep}(t,x=1,-\mu,\omega)\\
    & + 2(1-\eta_{\tr}(\omega)) f^\ep (t,x=1,\mu,\omega) \,, & \mu>0\\
    g^{\ep}(t,x=1+L_{\sub},\mu,\omega) =& 0\,, & \mu<0\\
    g^{\ep}(t=0,x,\mu,\omega) =& 0\,.
    \end{aligned}\right.
\end{equation}
%It is interesting to point out that it is the compatibility condition that enables the model reduction which intertwines the transducer with the substrate.
It is worth noting that the model reduction that intertwines the transducer with the substrate is made possible by the compatibility condition.

At the right boundary $x=1+L_{\sub}$, absorbing boundary condition is imposed for $g^{\ep}$. This is a safe assumption and was justified by the fact that the substrate is much thicker than the transducer, i.e., $L_{\sub}\gg 1$. See also~\cite{hua2017experimental}.

\section{Reconstruction formulation and main theorem}\label{sec:inverse_results}
We now describe the associated inverse problem. In particular, we discuss the experimental setting and, by defining a PDE solution map, set up an appropriate model for describing the inverse problem setup. Then in Section~\ref{subsec:stability_main_thm}, we examine the inverse problem stability and present the stability deterioration problem.

\subsection{Inverse problem setup}\label{subsec:inverse_setting}
In lab experiments, non-intrusive detection methods that do not disrupt the material under study are preferred. For an experiment at such a small scale, measurements can usually only be made on the surface. As presented in the introduction, laser beams are shot on the surface of the material that triggers heat propagation, and one could place temperature sensors on the surface of the material to detect the temperature fluctuations as a function of time.

% All other parameters in the equations are assumed to be known, and we seek to reconstruct the reflection coefficient of the transducer at the interface of the two bulk materials, namely $\eta_{\tr}$ in~\eqref{eq:forward_PTE_f_eps_couple}. 
Because the other coefficients at the interface, $\zeta_{\tr}, \eta_{\sub}, \zeta_{\sub}$, are determined by $\eta_{\tr}$ through~\eqref{eq:tr_sub_couple}~\eqref{eq:eta_1_compatible_condition}, it is sufficient to reconstruct the reflection coefficient of the transducer at the interface of the two bulk materials %, namely 
$\eta_{\tr}$ in~\eqref{eq:forward_PTE_f_eps_couple}. Henceforth, we will drop the subscript.

With this setup, we need to define a source-to-measurement operator as follows:
\begin{equation}\label{eq:forward_operator_M_eta}
\Lambda_{\eta}^\ep : \phi(t,\mu,\omega) \mapsto \Delta T(t,x=0)\,.
\end{equation}
This operator maps the input Dirichlet boundary condition $\phi$ to the temperature response $\Delta T(t)$ at the surface of the material as a function of time. Moreover, the measurement functional is defined as the projection of the operator~\eqref{eq:forward_operator_M_eta} onto some test function $\psi(t)$,
\begin{equation}\label{eq:forward_map_M_eta}
    \mathcal{M}^\ep(\eta): \phi(t,\mu,\omega)
    \mapsto
    \int\limits_{0}^{t_{\max}} \Lambda_{\eta}^\ep[\phi](t) \psi(t)\,\rd t = 
    \int\limits_{0}^{t_{\max}}\Delta T(t,x=0)\psi(t)\,\rd t\,,
\end{equation}
where $t_{\max}$ is the final time of the measurement that is determined by experimental setup. 
When all parameters are fixed, this map is uniquely determined by the reflection coefficient $\eta(\omega)$.

The inverse problem then amounts to:
\begin{equation}
\tag{Q}\label{Q:infer}
    \emph{inferring $\eta(\omega)$ using the measurement of $\mathcal{M}^\ep(\eta)$.}
\end{equation}

Numerically, this problem can be formulated as a PDE-constrained optimization. In experiments, one typically conducts a sequence of experiments with sources $\{\phi_i\}_{i=1}^I$. For each such experiment, we detect the temperature fluctuation using test functions $\{\psi_j\}_{j=1}^J$, and we can define the corresponding measurement functional by
\begin{equation}\label{eq:Meps_ij_def}
    \mathcal{M}^\ep_{i,j}(\eta):\phi_i(t,\mu,\omega)\mapsto\int\limits_{0}^{t_{\max}}\Delta T_i(t,x=0)\psi_j(t)\,\rd t = d_{i,j}\,.
\end{equation}
Then the optimization problem is written as:
\begin{equation}\label{eq:L_eta_PDE_optim}
    \eta(\omega)\in \argmin_{\eta \in\Omega}\,L^\ep(\eta),\quad L^{\ep}(\eta):=\frac{1}{IJ}\sum\limits_{i=1}^I\sum\limits_{j=1}^J|\mathcal{M}^\ep_{i,j}(\eta)-d_{i,j}|^2\,,
\end{equation}
where $\Omega$ is a feasible search space of $\eta$. This optimization problem is constrained by PTE such that $\mathcal{M}^\ep_{i,j}(\eta)$ is computed by the solution $f^{\ep}$ with the reflection coefficient $\eta$ through~\eqref{eq:forward_PTE_f_eps_couple}.

Theoretically, to answer~\eqref{Q:infer}, in~\cite{sun2022unique} the authors studied the uniqueness, using a singular decomposition technique. That is, to prepare the incoming data as a function singular in time, speed, and frequency. Then, along the characteristics, the singularity will be mostly preserved, and the remainder of the solution is proved negligible. The leading order (singular) part of the solution reveals the information about $\eta(\omega)$ through reflection at the right boundary. The paper concludes that for $\ep=1$, if $\sup_{\phi(\omega)\in \Phi}\left|\mathcal{M}^{\ep}(\eta_1)-\mathcal{M}^{\ep}(\eta_2)\right| = 0$ for $\Phi$ being the set of singular source function $\phi(\omega)$ concentrating around some $\omega_0\in(0,\infty)$, %\ql{in what sense}
then $\eta_1(\omega)=\eta_2(\omega)$ point-wise. %\ql{in what sense}.

The result does not contradict~\cite{hua2017experimental, Hua_Lindsay_general_Fourier_2019}%~\ql{the experimental result}.
, but does not answer the ill-conditioning observed in the same paper either. Through analysis, it will be seen that in the diffusive regime~\cite{sun2022unique}, %\ql{cite earlier paper},
although the inverse problem gives a unique reconstruction, the reconstruction is severely ill-posed in the sense of the lack of stability. In particular, we will show that for $\eta_1,\eta_2\in\Omega$, the difference $|\mathcal{M}^\ep(\eta_1)-\mathcal{M}^\ep(\eta_2)|$ can be arbitrarily small when $\|\eta_1-\eta_2\|$ is arbitrarily (with an exponential rate) large. We dedicate section~\ref{subsec:stability_main_thm} and section~\ref{sec:pf_stability_thm} for the theorem and the proof, and present our main assumptions here. 
\begin{assumption}\label{hyp1}
    \begin{enumerate}
    
    \item The phonon relaxation time is continuous and bounded from below
    \begin{equation}\label{eq:assump_rlx_bd}
        \tau(\omega)\ge \tau_{\min}>0,\quad\forall\omega \in (0,\infty)\,.
    \end{equation}
    \item The phonon group velocity is bounded and differentiable with
    \begin{equation}\label{eq:assump_vgroup_bd}
        0<\nu_{\min}\le \nu(\omega)\le \nu_{\max}<\infty,\quad\forall\omega \in (0,\infty)\,.
    \end{equation}

    \item The reflection index $\eta(\omega)$ belongs to the feasible set $\Omega:=C^0((0,\infty)\to [0,1])$.

    % \item The thickness of the substrate satisfies
    % \begin{equation}\label{eq:assump_Ls_thick}
    % L_s \geq \frac{\nu_{\max}}{2} + \frac{\nu_{\max}}{\mu_0 \nu(\omega_0)} + \frac{1}{2}\,.
    % \end{equation}
    
    \item There exist $p\in(1,3/2)$, and some $0<c_p<\infty$ such that
    \begin{equation}\label{eq:assump_bdd_moment}
        \int\limits_{0}^\infty\frac{C_\omega}{\tau^2(\omega)\nu(\omega)^{1/p}}\,\rd\omega\leq c_p<\infty\,.
    \end{equation}

    \item The source function $\phi \in L^1(\rd t\rd\mu\rd\omega)\cap L^p\left(0,t;L^p(C_{\omega}\rd\mu \rd\omega)\right)$ where $p$ is as in~\eqref{eq:assump_bdd_moment}, and there exists a constant $c_m>0$ such that $0\leq \phi \leq c_m C_{\omega}$. 
\end{enumerate}
\end{assumption}
Both the transducer and the substrate satisfy ~\eqref{eq:assump_rlx_bd}, ~\eqref{eq:assump_vgroup_bd} and~\eqref{eq:assump_bdd_moment} in assumption~\ref{hyp1}. The first two assumptions ensure that the relaxation time and the group velocity parameters in the PTE are finite, and hence the solution to the forward problem is well-defined. 
The last two assumptions require the physical quantities to have finite moments and the incoming source to be finite in some $L^p$-norm, which will contribute to the control of the solution by the incoming source.
Based on Assumptions~\ref{hyp1} and the feasible set of $\eta(\omega)$, we have that $f^\ep \in L^{p}(0,t_{\max};L^p(C_{\omega}\rd x\rd\mu \rd\omega))$ and continuous in~\eqref{eq:forward_PTE_f_eps_couple}~\cite{sun2022unique}.

\subsection{Stability of the inverse problem}\label{subsec:stability_main_thm}
One unique feature of this problem is that the hyperbolic component of the PDE has a straightforward characteristic line. So to reconstruct $\eta$ at a specific $\omega_0$, one can design a specific pair of source and measurement test functions to trigger the response at $\omega_0$. Hence, the stability of reconstruction of $\eta$ at $\omega_0$ depends solely on this one measurement. To be specific, let the
source function 
be singular in time $t$, direction cosine $\mu$, and phonon frequency $\omega$: %. Suppose that in~\eqref{eq:forward_PTE_f_eps_couple} the source is given by
\begin{equation}\label{eq:source_phi_t_mu_omega}
    \phi(t,\mu,\omega)=\frac{1}{\theta_\mu\theta_\omega\theta_t} \phi_t\left(\frac{t}{\theta_t}\right) \phi_\mu\left(\frac{\mu-\mu_0}{\theta_\mu}\right)\phi_\omega\left(\frac{\omega-\omega_0}{\theta_\omega}\right)\,,
\end{equation}
where $\phi_{\mu.\omega,t}(y)$ denote smooth functions supported in $|y|\leq1$, usually taken to be smooth cutoff functions. 
Correspondingly, the measurement test function takes the form of:
\begin{equation}\label{eqn:measurement_delta}
    \psi(t)=\psi_t\left(\frac{t-t_1}{\theta}\right), \quad t_1 = \frac{2\ep}{\mu_0\nu(\omega_0)}\,,
\end{equation}
with $\psi_t(t)$ is supported on $(-1,1)$. For these sources and test functions, we state the main stability estimate for one measurement:
\begin{theorem}\label{thm:stab_estim}
For any fixed $\omega_0$ and fixed $\varepsilon$ with the source in~\eqref{eq:source_phi_t_mu_omega} satisfying $\theta = \theta_{\omega}=\theta_{\mu}=\theta_t \leq \ep^{\alpha}$, with $\alpha>\frac{q+1}{q-3}$, and the test function in~\eqref{eqn:measurement_delta}, assuming Assumptions~\ref{hyp1} hold and that $t_{\max} > \frac{2\ep}{\nu(\omega_0)} + \theta$, there exists some constants $c_0, c_1>0$ independent of $\ep$, such that the measurement (defined in~\eqref{eq:forward_map_M_eta}) satisfies:
\begin{equation} \label{eq:measure_exp_estimate}
    |\mathcal{M}^\ep(\eta_1)-\mathcal{M}^\ep(\eta_2)|
    \leq c_1|\eta_1(\omega_0)-\eta_2(\omega_0)|e^{-c_0/\ep},\quad\text{as }\theta\to 0\,.
\end{equation}  
\end{theorem}
The inequality above shows that, even when one is given two distinct reflection coefficients $\eta_1$ and $\eta_2$ that can be far from each other, the sources and the measurements used in experiments cannot detect such differences: the difference between readings of the measurements is exponentially small when $\ep\to0$. Since the measurement data are incapable of extracting information from the reflection coefficient, this inequality essentially reveals the instability of the inversion.

This instability can be easily translated to the PDE-constrained optimization setting noting that the loss function $L^\ep[\cdot]$~\eqref{eq:L_eta_PDE_optim} is a finite sum of mismatches. Derived from this theorem, the following corollary is deduced to characterize the ill-conditioning of the corresponding optimization problem~\eqref{eq:L_eta_PDE_optim}.
\begin{corollary}\label{thm:cor_loss_eta_lipschitz}
Under the same assumption as in Theorem~\ref{thm:stab_estim} and denote the source function as
$\phi_i:=\frac{1}{\theta_\mu\theta_\omega\theta_t} \phi_t\left(\frac{t}{\theta_t}\right) \phi_\mu\left(\frac{\mu-\mu_i}{\theta_\mu}\right)\phi_\omega\left(\frac{\omega-\omega_i}{\theta_\omega}\right)$ 
with $\mu_i>0$ and $\psi_j:=\psi_t\left(\frac{t-t_j}{\theta}\right)$. There exist some constants $c_3=c_3(\eta_1, \eta_2, \{d_{i,j}, \phi_i, \psi_j\}_{i,j}), c_2=c_2(\{\phi_i\}_{i})$, such that the loss function of the PDE-constrained optimization~\eqref{eq:L_eta_PDE_optim} satisfies
\begin{equation}
|L^\ep(\eta_1)-L^\ep(\eta_2)| \leq c_3 \|\eta_1 -\eta_2\|_{L^\infty} e^{-c_2/\ep},\quad\text{as }\theta\to 0\,.
\end{equation}
\end{corollary}
Similar to the interpretation above, the Lipschitz constant of the loss function $L^\ep[\cdot]$ vanishes in the diffusion limit $\ep\to 0$. This suggests that the loss function is essentially flat and is incapable of differentiating different reflection coefficients in the diffusive limit.
\begin{proof}
To show this corollary, we will 
recall the definition of the measurement $\mathcal{M}^\ep_{i,j}$ associated with $(\phi_i,\psi_j)$ in~\eqref{eq:Meps_ij_def}, 
and apply the estimate~\eqref{eq:measure_exp_estimate} for every $\{i,j\}$. 
\begin{equation*}
\begin{aligned}
|L^\ep(\eta_1)-L^\ep(\eta_2)| = & \frac{1}{IJ}\sum_{i=1}^I \sum_{j=1}^J |\mathcal{M}^\ep_{i,j}(\eta_1) - \mathcal{M}^\ep_{i,j}(\eta_2)| |\mathcal{M}^\ep_{i,j}(\eta_1) + \mathcal{M}^\ep_{i,j}(\eta_2) - 2d_{ij}|\,,\\
\leq & \frac{1}{IJ} c_3 \sum_{i=1}^I \sum_{j=1}^J |\mathcal{M}^\ep_{i,j}(\eta_1) - \mathcal{M}^\ep_{i,j}(\eta_2)| \,,\\
\leq & \frac{1}{IJ} c_3 \sum_{i=1}^I \sum_{j=1}^J |\eta_1(\omega_i) - \eta_2(\omega_i)| e^{-c_{i}/\ep} \,, \\
\leq & c_3 \|\eta_1 - \eta_2\|_{L^\infty} e^{-c_2/\ep}\,,
\end{aligned}
\end{equation*}
where given fixed $\{\phi_i,\psi_j\}$, $c_3$ is some uniform bound of the measurement operator and the data, 
$c_i$ is a constant depending on $\{\mu_i,\omega_i, \nu(\omega_i)\}$,
and $c_2$ is determined by $\max_i\{\mu_i\}, \nu_{\max}, \tau_{\max}$.
\end{proof}

The proof of Theorem~\ref{thm:stab_estim} is built upon the singular decomposition technique. This method is to design a mechanism to separate the ballistic and scattering parts of the solution, and trace the effect of each separately. More specifically, let $f^\ep$ solve the original phonon transport equation~\eqref{eq:forward_PTE_f_eps_couple} and define $f^\ep=:f_0^\ep+f_1^\ep$, with $f_0^\ep$ solving the ballistic part:
\begin{equation}\label{eqn:g0}
\left\{\begin{aligned}
\ep\partial_t f_0^\ep+\mu \nu \partial_x f_0^\ep &=-\dfrac{1}{\ep\tau}f_0^\ep\,, \\
f_0^\ep(t,x=0,\mu,\omega) &= {\phi}(t,\mu,\omega)\,, & \mu>0\,,\\
f_0^\ep(t,x=1,\mu,\omega) &= \eta(\omega)f_0^\ep(t,x=1,-\mu,\omega)\,, & \mu<0\,, \\
f_0^\ep(t=0,x,\mu,\omega) &=0\,.
\end{aligned}\right.
\end{equation}
Then $f_1^\ep$ that keeps track of the remainders solves the scattering part:
\begin{equation}\label{eqn:g1}
\left\{\begin{aligned}
\ep\partial_t f_1^\ep+\mu \nu \partial_xf_1^\ep & = \dfrac{1}{\ep\tau}\left(-f_1^\ep + \frac{C_\omega}{C_\tau}\langle f^\ep/\tau\rangle\right)\,,\\
f_1^\ep(t,x=0,\mu,\omega) &= 0 \,, &\mu>0\,,\\
f_1^\ep(t,x=1,\mu,\omega) &=\eta(\omega) f_1^\ep(t,x=1,-\mu,\omega) 
+ \zeta(\omega) g^{\ep}(t,x=1,\mu,\omega)\,,& \mu<0\,, \\
f_1^\ep(t=0,x,\mu,\omega)& =0\,.
\end{aligned}\right.
\end{equation}
Since the incoming data $\phi$ is very singular, the ballistic component of the solution $f_0^\ep$ will keep track of such singularity, and will contribute largely to the measurement. In comparison, $f_1^\ep$ is spread out, and its contribution can be small. The main component of the proof is to demonstrate such a separation of contribution.
% This decomposition allows us to decompose the measurement operator and obtain the following estimates. Along the characteristics, we will show that the singularity in frequency $\omega_0$ of the incoming sources is preserved by $f_0$, and the measurement is mainly from $f_0$, which contains information about the transmission index by reflection at the boundary. 
\begin{theorem}\label{them:measurem_decomp}
Assuming Assumptions~\ref{hyp1} hold, using the data prepared in~\eqref{eq:source_phi_t_mu_omega} and~\eqref{eqn:measurement_delta},
the measurement~\eqref{eq:forward_map_M_eta} admits the decomposition
\begin{equation}\label{eq:Measure_decompose_12}
    \mathcal{M}^\ep (\eta) = \mathcal{M}_0^\ep (\eta) + \mathcal{M}_1^\ep (\eta)\,,
\end{equation}
where each term is defined as follows, %we recall the definitions of temperature~\eqref{eq:DeltaT_def} and the measurement
\begin{equation}\label{eq:M^eps_measurement}
\mathcal{M}_j^\ep(\eta)=\frac{1}{C_\tau}\int\limits_{0}^{t_{\max}}\psi(t)\int\limits_{-1}^1\int\limits_0^\infty f^\ep_j(t,x=0,\mu,\omega)/\tau(\omega)\, \rd\omega\rd\mu\rd t,\quad j=0,1\,.
\end{equation}
Moreover, there exist positive constants $c_1, c_4, c_5$ such that
\begin{equation}\label{eq:M0_converge_theta}
\lim_{\theta\to 0} \mathcal{M}_0^\ep(\eta) = c_1\eta(\omega_0)e^{-\frac{2}{\mu_0\nu_0\tau_0\ep}} + c_4,\quad
\text{for }\theta_\omega=\theta^{\alpha_1}, \theta_\mu=\theta^{\alpha_2}, \alpha_{1,2}\geq 1, \theta_t=\theta \,, 
\end{equation}
where $c_1, c_4$ depend on the choice of $\alpha_{1,2}$, but are independent of $\eta$, and
\begin{equation}\label{eq:M1_eps_estimate}
\mathcal{M}_1^\ep(\eta)\le c_5 \frac{\theta}{(\theta_\mu\theta_\omega\theta_t)^{1/q}\ep^{1+1/q}}\,,
\end{equation}
where $q$ is the H\"{o}lder conjugate of $p$ in assumption~\ref{hyp1} that $\frac{1}{q}+\frac{1}{p}=1$.
\end{theorem}

The stability estimate in Theorem~\ref{thm:stab_estim} is a direct consequence of Theorem~\ref{them:measurem_decomp}. With assumptions in Theorem~\ref{thm:stab_estim}, 
% In the case where $\theta = \theta_{\omega}=\theta_{\mu}=\theta_t \leq \eps^{\alpha}$, with $\alpha>\frac{q+1}{q-3}$, 
we have both $\mathcal{M}_1^\ep(\eta_1), \mathcal{M}_1^\ep(\eta_2)\to 0$, as $\theta\to 0$, while the difference 
between $\mathcal{M}_0^\ep(\eta_1), \mathcal{M}_0^\ep(\eta_2)$ is estimated by
$|\mathcal{M}_0^\ep(\eta_1)-\mathcal{M}_0^\ep(\eta_2)|\leq c_1|\eta_1(\omega_0)-\eta_2(\omega_0)| e^{-\frac{2}{\mu_0 \nu_0 \tau_0 \ep}}$, and Theorem~\ref{thm:stab_estim} is concluded by the triangle inequality.

\section{Proof of the stability estimate}\label{sec:pf_stability_thm}

We will prove Theorem~\ref{them:measurem_decomp} in this section. More specifically, we will trace the contribution towards the measurement of the singular component $f_0^\ep$ and the scattering component $f_1^\ep$. In the $\theta\to0$ limit, we will show that the scattering component's contribution vanishes, and that the singular component's contribution loses sensitivity. The estimates of these two components are shown respectively in section~\ref{subsec:measure_f0} and~\ref{subsec:measure_f1}.

Before we proceed, we first construct the solution to~\eqref{eqn:g0} and~\eqref{eqn:g1} in the following proposition.
\begin{proposition}\label{prop:charac}
The solution to~\eqref{eq:forward_PTE_f_eps_couple} admits the decomposition $f^\ep=f^\ep_0+f^\ep_1$ such that at $x=0$ and $\mu<0$, the solution reads
\begin{itemize}
\item For $t\geq 0$,
\begin{equation}\label{eqn:g0_character_neg}
f_0^\ep(t,x=0,\mu<0,\omega) = \eta(\omega)\phi\Big(t-\frac{2\ep}{|\mu|\nu}, |\mu|,\omega\Big)e^{-\frac{2}{\ep|\mu|\nu\tau}}\mathbf{1}_{t>\frac{2\ep}{|\mu|\nu}}\,.
\end{equation}

\item For $0\leq t \leq \frac{\ep}{|\mu|\nu}$,
\begin{equation}
f_1^\ep(t,x=0,\mu<0,\omega)=\frac{1}{\ep^2\tau} \frac{C_\omega}{C_\tau} \int\limits_0^{t}\langle f^\ep/\tau\rangle \left(t-s,-\frac{\mu\nu s}{\ep}\right) e^{-\frac{s}{\ep^2\tau}}\,\rd s
=: f_{1,0}^\ep\,.
\end{equation}

\item For $\frac{\ep}{|\mu|\nu}<t \leq \frac{2\ep}{|\mu|\nu}$,
\begin{equation}\label{eq:f1_123}
\begin{aligned}
& f_1^\ep(t,x=0,\mu<0,\omega)\\
=&\, \frac{\eta(\omega)}{\ep^2\tau} \frac{C_\omega}{C_\tau} e^{-\frac{1}{\ep|\mu|\nu \tau}} \int\limits_0^{t-\frac{\ep}{|\mu|\nu}} \langle f^\ep/\tau \rangle\left(t-s-\frac{\ep }{|\mu|\nu}, 1+\frac{\mu\nu s}{\ep}\right) e^{-\frac{s}{\ep^2\tau}}\, \rd s\\
& + \zeta(\omega) g^{\ep} \left(t-\frac{\ep}{|\mu| \nu},1,\mu,\omega\right) e^{-\frac{1}{\ep^2\tau} \frac{\ep}{|\mu|\nu}}\\
& +\frac{1}{\ep^2 \tau} \frac{C_\omega}{C_\tau} \int\limits_{0}^{\frac{\ep}{|\mu|\nu}}\langle f^\ep/\tau \rangle \left(t-s, x-\frac{\mu\nu s}{\ep}\right)e^{-\frac{s}{\ep^2\tau}}\,\rd s\,,\\
=: & f_{1,1}^\ep + f_{1,2}^\ep + f_{1,3}^\ep \,.
\end{aligned}
\end{equation}

\item For $t > \frac{2\ep}{|\mu|\nu}$,
\begin{equation}\label{eq:f1_456}
\begin{aligned}
f_1^\ep(t,x=0,\mu<0,\omega)=
& \frac{\eta(\omega)}{\ep^2\tau} \frac{C_\omega}{C_\tau} e^{-\frac{(1-x)}{\ep|\mu|\nu\tau}} \int\limits_0^{\frac{\ep}{|\mu|\nu}}\langle f^\ep/\tau \rangle \left(t-s-\frac{\ep}{|\mu|\nu}, 1+\frac{\mu\nu s}{\ep}\right) e^{-\frac{s}{\ep^2\tau}}\,\rd s\\
& + \zeta(\omega) g^{\ep} \left(t-\frac{\ep}{|\mu| \nu},1,\mu,\omega\right) e^{-\frac{1}{\ep^2\tau} \frac{\ep}{|\mu|\nu}}\\
& + \frac{1}{\ep^2\tau} \frac{C_\omega}{C_\tau} \int\limits_{0}^{\frac{\ep}{|\mu|\nu}}\langle f^\ep/\tau \rangle \left(t-s, x-\frac{\mu\nu s}{\ep} \right)e^{-\frac{s}{\ep^2\tau}}\,\rd s\,,\\
=: & f_{1,4}^\ep + f_{1,5}^\ep + f_{1,6}^\ep\,.
\end{aligned}
\end{equation}
\end{itemize}
\end{proposition}
This proposition is a direct application of the method of characteristics. We omit the proof here for brevity and detail the calculations in Appendix~\ref{subsec:charac}.
These explicit formulas enable our estimates for the two parts of temperature measurement, $\mathcal{M}_0^\ep$ and $\mathcal{M}_1^\ep$.

\subsection{Estimate of \texorpdfstring{$\mathcal{M}^\ep_0$}{Meps0}} \label{subsec:measure_f0}
We provide the estimate for the singular component $\mathcal{M}_0^\ep$ in this section. Note that $f^\ep_0$ has $\mu>0$ and $\mu<0$ two components:
\[
\mathcal{M}_0^\ep 
=: \mathcal{M}_{0,\mu<0}^\ep + \mathcal{M}_{0,\mu>0}^\ep\,.
\]
To control $\mathcal{M}_{0,\mu>0}^\ep$, we plug in the incoming boundary of $f_0^\ep$~\eqref{eqn:g0} for:
\begin{equation*}
\begin{aligned}
\mathcal{M}_{0,\mu>0}^\ep & = \frac{1}{C_\tau} \lim_{\theta,\theta_{\mu},\theta_{\omega}\to0} \int\limits_{0}^{\infty} \int\limits_{0}^{\infty} \int\limits_{-1}^1 \frac{\phi_\mu(\mu')\phi_\omega(\omega')\phi_t(t')}{\tau(\omega'\theta_\omega+\omega_0)}  \psi_t\left(t'-\frac{2\ep}{\theta \mu_0 \nu(\omega_0)}\right) \,\rd\mu'\rd\omega' \rd t'\\
& =: c_4\,,
\end{aligned}
\end{equation*}
which is independent of $\eta(\omega)$ and gives the constant $c_4$ in Theorem~\ref{them:measurem_decomp}.

To estimate $\mathcal{M}_{0,\mu<0}^\ep$, we plug in~\eqref{eqn:g0_character_neg} for:
\begin{equation}
\begin{aligned}
\mathcal{M}^{\ep}_{0,\mu<0} := \frac{1}{\theta_\mu\theta_\omega\theta_t C_\tau}
\int\limits_{-1}^0 \int\limits_{0}^{\infty}
\int\limits_{\frac{2\ep}{\mu\nu}}^{t_{\max}}
\frac{\eta(\omega)}{\tau(\omega)} 
& \phi_\mu\left(\frac{|\mu|-\mu_0}{\theta_\mu}\right)
\phi_\omega\left(\frac{\omega-\omega_0}{\theta_\omega}\right)
\phi_t\left(\frac{t}{\theta_t}-\frac{2\ep}{|\mu|\nu \theta_t}\right)\\
& \times e^{-\frac{2}{|\mu|\nu \tau \ep}}
\psi_t\left(\frac{t-t_1}{\theta}\right)\,\rd t\rd\omega \rd\mu\,.  
\end{aligned}
\end{equation}
Using the change of variables
%\begin{equation}
$\mu':=\frac{|\mu|-\mu_0}{\theta_\mu} = \frac{-\mu-\mu_0}{\theta_\mu}\,, 
\omega':= \frac{\omega-\omega_0}{\theta_\omega}\,,
t' := \frac{1}{\theta_t}\left(t-\frac{2\ep}{|\mu|\nu}\right)\,,$
%{\mu}\to\frac{\mu-\mu_0}{\theta_\mu}, \quad 
%{\omega}\to\frac{\omega-\omega_0}{\theta_\omega},\quad 
%{t}\to\frac{1}{\theta_t}\left(t-\frac{2\eps}{\mu\nu}\right)\,,
%\end{equation}
we can rewrite the above integral as
\begin{equation}
\mathcal{M}^{\ep}_{0,\mu<0} = \frac{1}{C_\tau}
\int\limits_{\frac{-\mu_0}{\theta_\mu}}^{\frac{1-\mu_0}{\theta_\mu}} \phi_\mu(\mu')  
\int\limits_{-\frac{\omega_0}{\theta_\omega}}^{\infty} \frac{\eta(\omega_0+\theta_\omega\omega')}{\tau(\omega_0+\theta_\omega\omega')} \phi_\omega(\omega') 
e^{-\frac{2}{\ep \mu_1 \nu_1 \tau_1}}
\int\limits_{0}^{t_{\text{up}}}
\phi_t(t')
\psi(t_{\text{test}})\,\rd t' \rd\omega' \rd\mu'\,,
\end{equation}
where 
%\begin{equation} 
$\mu_1 = \theta_{\mu}\mu' + \mu_0\,,
\nu_1 = \nu(\theta_{\omega}\omega' + \omega_0)\,, 
\tau_1 = \tau(\theta_{\omega}\omega' + \omega_0)\,,$ 
%\end{equation}
and 
%\begin{equation}
$t_{\text{up}}:=\frac{1}{\theta_t}\left(t_{\max} - \frac{2\ep}{\mu_1 \nu_1}\right)$, 
$t_{\text{test}}:=\frac{\theta_t t'}{\theta} + \frac{1}{\theta}\frac{2\ep}{\mu_1\nu_1} - \frac{t_1}{\theta}.$
%\end{equation}
% \begin{equation}
% \begin{aligned}
% \int\limits_{\mu=-\frac{\mu_0}{\theta_\mu}}^{1-\frac{\mu_0}{\theta_\mu}}
% \int\limits_{\omega=-\frac{\omega_0}{\theta_\omega}}^{\infty}
% \int\limits_{t=0}^{t_{\mathrm{up}}}
% \frac{\eta(\omega_0+\theta_\omega{\omega})}{\tau(\omega_0+\theta_\omega{\omega})}\phi_\mu({\mu})\phi_\omega({\omega})
% \phi_t(t)\psi_t(t')
% e^{-\frac{2}{\eps(\mu')\nu(\omega')\tau(\omega')}}\,\rd\mu\rd\omega  \rd t\,,
% \end{aligned}
% \end{equation}
Recall that $t_1=\frac{2\ep}{\mu_0\nu(\omega_0)}$,
and for the specific choice that
$\theta_t=\theta$, %$\theta_\omega=\theta_\mu$
by Taylor expansion,
\begin{equation}
\begin{aligned}
t_{\text{test}} =&\, t' - \frac{t_1}{\theta} + \frac{2\ep}{\theta(\mu_0 + \theta_\mu \mu')\nu(\omega_0 + \theta_\omega \omega')}\,, \\
=&\, t' + \frac{1}{\theta}\left(\frac{2\ep}{\mu_0\nu(\omega_0)}-t_1\right) - \frac{2\ep}{[\mu_0\nu(\omega_0)]^2} \left(\frac{\theta_\mu}{\theta}\mu'\nu(\omega_0) + \frac{\theta_\omega}{\theta} \omega'\mu_0\nu'(\omega_0)\right) \\
&\, + O\left(\frac{\theta_\mu^2 + \theta_\omega^2}{\theta}\right)\,. 
\end{aligned}
\end{equation}
% From here, we see that we need $t_1=\frac{2\eps}{\mu_0\nu(\omega_0)}$, $\theta_\mu,\theta_\omega=\theta^{\alpha_j}$ and $\theta_t=\theta$ for a non-trivial value for $\mathcal{M}_0$ 
% due to the compact support of $\psi_t$.  This gives 
% \begin{equation}
% \begin{aligned}
% \lim_{\theta_\mu, \theta_\omega\to 0}\frac{2\eps}{(\mu_0+\theta_\mu\mu)\nu(\omega_0+\theta_\omega\omega)}-t_1&=\lim_{\theta_\mu\to 0}\frac{2\eps}{(\mu_0+\theta_\mu\mu)\nu(\omega_0+\theta_\omega\omega)}-\frac{2\eps}{\mu_0\nu(\omega_0)}\\
%     &=-\frac{2\eps(\nu(\omega_0){\mu}+\mu_0\nu'(\omega_0){\omega})}{(\mu_0\nu(\omega_0))^2}:=\eps t^\ast(\mu,\mu_0,\omega_0)\,.
% \end{aligned}
% \end{equation}
Moreover, 
\begin{equation}
    \lim_{\theta_{\mu}, \theta_\omega\to 0} \frac{\eta(\omega_0+\theta_\omega \omega')}{\tau(\omega_0+\theta_\omega\omega')} e^{-\frac{2}{\ep \mu_1 \nu_1 \tau_1}}=\frac{\eta(\omega_0)}{\tau(\omega_0)}e^{-\frac{2}{\ep\mu_0\nu(\omega_0)\tau(\omega_0)}}\,.
\end{equation}
Hence, by using the continuity of functions (Assumptions~\ref{hyp1}) and the dominated convergence theorem (DCT), 
\begin{equation}
    \lim_{\theta\to 0}\, \mathcal{M}^{\ep}_{0,\mu<0}= c_1\eta(\omega_0)e^{-\frac{2}{\mu_0\nu_0\tau_0\ep}}\,,
\end{equation}
where $\tau_0:=\tau(\omega_0)$.
For $\theta_{\omega}=\theta^{\alpha_1}, \theta_\mu=\theta^{\alpha_2}$ as in~\eqref{eq:M0_converge_theta}, the constant $c_1$ is estimated according to the following cases:
\begin{itemize}
    \item If $\alpha_1, \alpha_2>1$, we have $\lim\limits_{\theta_{\mu},\theta_{\omega}\to0} \, t_{\text{test}} = t'$, and the constant $c_1$ is given by
    \begin{equation}
    \begin{aligned}
c_1 &= \frac{1}{\tau(\omega_0)}\frac{1}{C_\tau} \int\limits_{0}^{\infty} \int\limits_{0}^{\infty} \int\limits_{-1}^{1} \phi_\mu({\mu'})\phi_\omega({\omega'}) \phi_t(t')\psi_t(t') \,\rd\mu'\rd\omega' \rd t'\,.
\end{aligned}
\end{equation}
    \item If $\alpha_1= \alpha_2=1$, i.e., $\theta_\mu=\theta_{\omega}=\theta$ as in Theorem~\ref{thm:stab_estim}, then $\lim\limits_{\theta_{\mu},\theta_{\omega}\to0}\, t_{\text{test}} = t'-\frac{2\ep}{[\mu_0\nu(\omega_0)]^2} \left(\mu'\nu(\omega_0) + \omega'\mu_0\nu'(\omega_0)\right)$, and
    \begin{equation}
    \begin{aligned}
    c_1 =&\, \frac{1}{\tau(\omega_0)}\frac{1}{C_\tau} \int\limits_{0}^{\infty} \int\limits_{0}^{\infty} \int\limits_{-1}^1 \phi_\mu({\mu'})\phi_\omega({\omega'}) \phi_t(t')\\
    & \,\times \psi_t \left(t'-\frac{2\ep}{[\mu_0\nu(\omega_0)]^2} \left(\mu'\nu(\omega_0) + \omega'\mu_0\nu'(\omega_0)\right)\right) \,\rd\mu'\rd\omega' \rd t'\,.
    \end{aligned}
    \end{equation}
    \item If $\alpha_1=1, \alpha_2>1$, then $\lim\limits_{\theta_{\mu},\theta_{\omega}\to0}\, t_{\text{test}} = t'-\frac{2\ep}{[\mu_0\nu(\omega_0)]^2} \omega'\mu_0\nu'(\omega_0)$, and
    \begin{equation}
    \begin{aligned}
    c_1 &= \frac{1}{\tau(\omega_0)}\frac{1}{C_\tau} \int\limits_{0}^{\infty} \int\limits_{0}^{\infty} \int\limits_{-1}^1 \phi_\mu({\mu'})\phi_\omega({\omega'}) \phi_t(t')\psi_t \left(t'-\frac{2\ep}{[\mu_0\nu(\omega_0)]^2} \omega'\mu_0\nu'(\omega_0)\right) \,\rd\mu'\rd\omega' \rd t'\,.
    \end{aligned}
    \end{equation}
    \item If $\alpha_1>1, \alpha_2=1$, the estimate is similar to the previous case with the limit of $t_{\text{test}}$ replaced by $\lim\limits_{\theta_{\mu},\theta_{\omega}\to0}\, t_{\text{test}} = t'-\frac{2\ep}{[\mu_0\nu(\omega_0)]^2}\mu'\nu(\omega_0)$.
\end{itemize}

This concludes the estimate in~\eqref{eq:M0_converge_theta}.
% And the constant $C_3$ in Theorem~\ref{them:measurem_decomp} is given by the incoming data:
% \begin{equation*}
% C_3 = \frac{1}{C_\tau} \lim_{\theta,\theta_{\mu},\theta_{\omega}\to0} \int\limits_{0}^{\infty} \int\limits_{0}^{\infty} \int\limits_{-1}^1 \frac{1}{\tau(\omega'\theta_\omega+\omega_0)} \phi_\mu(\mu')\phi_\omega(\omega')\phi_t(t') \psi_t\left(t'-\frac{2\eps}{\theta \mu_0 \nu(\omega_0)}\right) \,\rd\mu'\rd\omega' \rd t'\,,
% \end{equation*}
% which does not depend on $\eta(\omega)$.

\subsection{Estimate of \texorpdfstring{$\mathcal{M}^\ep_1$}{Meps1}} \label{subsec:measure_f1} 

Now, we are to show that the estimate of measurement from $f_1^\ep$~\eqref{eq:M1_eps_estimate}. Recall that $f_1^\ep$ satisfies zero incoming boundary condition~\eqref{eqn:g1} and the decomposition of $f_1^\ep$ in~\eqref{eq:f1_123} and~\eqref{eq:f1_456}, we have:
\[
\mathcal{M}_1^{\ep} = \mathcal{M}^{\ep}_{1,\mu<0} =: \sum_{i=0}^6 \mathcal{M}^{\ep}_{1,i}\,. 
\]
The rest of the section is dedicated to showing that each term in this decomposition provides a negligible contribution.

We should first emphasize that, according to this decomposition, the components $f^\ep_{1,0}$, $f^\ep_{1,1}$, $f^\ep_{1,3}$, $f^\ep_{1,4}$ and $f^\ep_{1,6}$ see the contribution of the scattering effects brought by $\langle\cdot\rangle$, while the components $f^\ep_{1,2}$ and  $f^\ep_{1,5}$ see the contribution of $g^{\ep}$.

Therefore, in what follows, we need to separately bound the scattering and $g^\ep$ terms.
\begin{lemma}[Proposition 3.1 in~\cite{sun2022unique}] \label{lem:prop_31_Lp_norm}
Let $\phi$ be an incoming data in~\eqref{eq:forward_PTE_f_eps_couple}, and $(f^\ep, g^\ep)$ be a solution to the system~\eqref{eq:forward_PTE_f_eps_couple} and~\eqref{eq:forward_PTE_g_eps_couple}, then for any $p\geq 1$, there exists positive constant $c_6=(1+c)^{1/p}$ such that
\begin{equation}\label{eq:f_g_Lp_estimate}
\|f^\ep(t,\cdot)\|_{L^p(C_{\omega}^{1-p}\,\rd\omega\rd\mu\rd x)} + \|g^\ep(t,\cdot)\|_{L^p(C_{\omega}^{1-p}\,\rd\omega\rd\mu\rd x)} 
\leq c_6 \|\phi\|_{L^p(0,t;L^p(C_{\omega}^{1-p}\,\rd\omega\rd\mu))}\,,
\end{equation}
here $\|\cdot\|_{L^p(C_{\omega}^{1-p}\,\rd\omega\rd\mu\rd x)}$ denotes the $L^p$-norm in $(x,\mu,\omega)$ with weight $C_{\omega}^{1-p}$.
\end{lemma}

The control of $\langle f^\ep/\tau \rangle$ is summarized in the following proposition:

\begin{proposition}\label{prop:4.3}
There exists some positive constant $c_7=c_7({\phi,\psi,p,\tau_{\min}, C_\tau})$ such that
\begin{equation} \label{eq:estimate_temp_feps}
\begin{aligned}
& \frac{1}{C_\tau^2} \frac{1}{\ep^2} 
\int\limits_{0}^{t_{\max}} \int\limits_{\mu<0} \int\limits_{\omega} \int\limits_{0}^t
\psi_t\left(\frac{t-t_1}{\theta}\right) 
\frac{C_\omega}{\tau^2(\omega)} \langle |f^\ep|/\tau \rangle\left(t-s,-\frac{\mu\nu s}{\ep}\right) \mathbf{1}_{\frac{|\mu|\nu s}{\ep}<1} \,\rd s\rd \omega \rd \mu \rd t \\
\le &\, c_7\theta(\theta_\mu\theta_\omega\theta_t)^{-\frac{1}{q}} \ep^{-\left(1+\frac{1}{q}\right)}\,.
\end{aligned}
\end{equation}
% \begin{equation}
% \int\limits_{0}^{t_{\max}} \int\limits_{\mu<0} \int\limits_{\omega} \int\limits_{0}^t\psi_t\Big(\frac{t-t_1}{\theta}\Big)\frac{C_\omega}{\eps^2\tau^2(\omega)}\langle |f^\eps|\rangle\left(t-s,-\frac{\mu\nu s}{\eps}\right) \mathbf{1}_{\frac{|\mu|\nu s}{\eps}<1} \,\rd s\rd \omega \rd \mu \rd t
% \le C_{\phi,\psi,p} \theta(\theta_\mu\theta_\omega\theta_t)^{-1/q} \eps^{-\big(1+\frac{1}{q}\big)} \,.
% \end{equation}
Here, $q$ is the H\"{o}lder conjugate of $p$ as in assumption~\ref{hyp1}.
\end{proposition}
The proof of this proposition is in Appendix~\ref{apx:scatter_theta_bound}. The estimate in~\eqref{eq:estimate_temp_feps} shows that the scattering part is controlled by the singular sources, indicated by $\theta$'s. Under the diffusive scaling, the rescaling of both time and space introduces an explicit dependence on the Knudsen number $\varepsilon$ in the upper bound.

The direct application of this proposition provides the estimate for $\mathcal{M}^{\ep}_{1, i}$ for $i=0,1,3,4,6$. Indeed,
\begin{itemize}
\item[--] $\mathcal{M}^{\ep}_{{1,0}}$:
\begin{equation}
\begin{aligned}
\mathcal{M}^{\ep}_{{1,0}} :=& 
\frac{1}{C_\tau} \int\limits_{t=0}^{t_{\max}} \psi_t\left(\frac{t-t_1}{\theta}\right) \int\limits_{\mu<0,\omega} \frac{1}{\ep^2\tau^2} \frac{C_\omega}{C_\tau} \\
& \times \int\limits_{s=0}^t \langle f^\ep/\tau \rangle \left(t-s,-\frac{\mu\nu s}{\ep}\right)e^{-\frac{s}{\ep^2\tau}}\, \rd s\,\mathbf{1}_{t\in(0,\frac{\ep}{|\mu|\nu})}\,\rd\mu\rd\omega\rd t \\
\le &\, c_7 \theta (\theta_\mu \theta_\omega \theta_t)^{-\frac{1}{q}} \ep^{-\left(1+\frac{1}{q}\right)}\,.
\end{aligned}
\end{equation}

\item[--] $\mathcal{M}^{\ep}_{{1,1}}$:
\begin{equation}
\begin{aligned}
    \mathcal{M}^{\ep}_{1,1} & := 
    \frac{1}{C_\tau^2} \int\limits_{t=0}^{t_{\max}}\psi_t\left(\frac{t-t_1}{\theta}\right)\int\limits_{\omega} \int\limits_{-1}^0 
    \int\limits_{s=0}^{t+\frac{\ep}{\mu\nu}}
    \eta(\omega) e^{-\frac{1}{|\mu|\nu\ep\tau}} e^{-\frac{s}{\ep^2\tau}}
    \frac{C_\omega}{\ep^2\tau^2} \\
    & \hspace{1.2in} \times \langle f^\ep/\tau \rangle \left(t-s+\frac{\ep}{\mu\nu},1+\frac{\mu\nu}{\ep}s\right)\, \rd s\, \mathbf{1}_{t\in(\frac{\ep}{|\mu|\nu},\frac{2\ep}{|\mu|\nu})}\, \rd\mu\rd\omega\rd t\,, \\
    &\le \frac{1}{C_\tau^2} \int\limits_{t=0}^{t_{\max}} \psi_t\left(\frac{t-t_1}{\theta}\right)
    \int\limits_{\omega} \int\limits_{-1}^0 \int\limits_{s=0}^{t+\frac{\ep}{\mu\nu}} \frac{C_\omega}{\ep^2\tau}\langle |f^\ep|/\tau \rangle \left(t -s +\frac{\ep}{\mu\nu}, 1+\frac{\mu\nu}{\ep}s \right)\,\rd s\\
    & \hspace{3.2in} \times 
    \mathbf{1}_{t\in(\frac{\ep}{|\mu|\nu},\frac{2\ep}{|\mu|\nu})} \, \rd\mu\rd\omega\rd t\,,
\end{aligned}
\end{equation}
where we used $\eta(\omega)\leq 1$, and $e^{-\frac{1}{\ep|\mu|\nu\tau}}, e^{-\frac{s}{\ep^2\tau}}\leq 1$.
Using the change of variables $s':=s-\frac{\ep}{\mu\nu}$, we obtain:
\begin{equation}
\begin{aligned}
\mathcal{M}^{\ep}_{1,1} 
& \leq  
    \frac{1}{C_\tau^2} \int\limits_{t=0}^{t_{\max}} \psi_t\left(\frac{t-t_1}{\theta}\right)
    \int\limits_{\omega} \int\limits_{-1}^0 
    \int\limits_{\frac{\ep}{|\mu|\nu}}^{t}
    \frac{C_\omega}{\ep^2\tau^2} \langle |f^\ep|/\tau \rangle \left(t-s', 2+\frac{\mu\nu}{\ep}s'\right)\, \rd s'\\
    & \hspace{3in} \times \mathbf{1}_{t\in(\frac{\ep}{|\mu|\nu},\frac{2\ep}{|\mu|\nu})}\, \rd\mu\rd\omega\rd t\,,
\end{aligned}
\end{equation}
Now let $\mu':=-\mu - \frac{2\ep}{\nu s'}=-\mu \left(1-\frac{2\ep}{|\mu|\nu s'}\right)$, 
then $2+\frac{\mu\nu}{\ep}s'=-\frac{\mu' \nu s'}{\ep}$. 
Moreover, since $\frac{\ep}{|\mu|\nu}<s'<t<\frac{2\ep}{|\mu|\nu}$,
we have $-1<1-\frac{2\ep}{|\mu|\nu s'}<0$, $\mu'\in(-1,0)$, and $\frac{|\mu'|\nu s'}{\ep}<1$. Therefore,
\begin{equation}
\begin{aligned}
\mathcal{M}^{\ep}_{1,1} & \leq  
    \frac{1}{C_\tau^2} \int\limits_{t=0}^{t_{\max}} \psi_t\left(\frac{t-t_1}{\theta}\right)
    \int\limits_{\omega} \int\limits_{-1}^0 
    \int\limits_{0}^{t}
    \frac{C_\omega}{\ep^2\tau^2} \langle |f^\ep|/\tau \rangle \left(t-s', -\frac{\mu'\nu s'}{\ep}\right) \\
    & \hspace{3in} \times \mathbf{1}_{\frac{|\mu'|\nu s'}{\ep}<1} \, \rd s'\, \rd\mu'\rd\omega\rd t \,,\\
    & \le c_7 \theta (\theta_\mu \theta_\omega \theta_t)^{-\frac{1}{q}} \ep^{-\left(1+\frac{1}{q}\right)}\,.
\end{aligned}
\end{equation}
\item[--] $\mathcal{M}^{\ep}_{{1,3}}$:
\begin{equation}
    \begin{aligned}
        \mathcal{M}^{\ep}_{1,3} &:= \frac{1}{C_\tau} \int\limits_{t=0}^{t_{\max}} \psi_t\left(\frac{t-t_1}{\theta}\right) 
        \int\limits_{\mu<0,\omega} \int\limits_{s=0}^{\frac{\ep}{|\mu|\nu}} 
        \frac{1}{\ep^2\tau^2} \frac{C_\omega}{C_\tau}\langle f^\ep/\tau \rangle \left(t-s,-\frac{\mu\nu}{\ep}s \right)e^{-\frac{s}{\ep^2\tau}} \,\rd s \\
        & \hspace{3in} \times\mathbf{1}_{t\in(\frac{\ep}{|\mu|\nu},\frac{2\ep}{|\mu|\nu})}\,\rd\mu \rd\omega \rd t\,,\\
        &\le \frac{1}{C_\tau^2} \int\limits_{t=0}^{t_{\max}} \psi_t\left(\frac{t-t_1}{\theta}\right) \int\limits_{\mu<0,\omega} \int\limits_{s=0}^{t} \frac{C_\omega}{\ep^2\tau^2}\langle |f^\ep|/\tau\rangle \left(t-s,-\frac{\mu\nu}{\ep}s\right)\mathbf{1}_{\frac{|\mu|\nu s}{\ep}<1}\, \rd s \rd\mu \rd\omega \rd t\,,\\
        & \le c_7 \theta (\theta_\mu \theta_\omega \theta_t)^{-\frac{1}{q}} \ep^{-\left(1+\frac{1}{q}\right)}\,.
    \end{aligned}
\end{equation}
\end{itemize}

Expanding $\mathcal{M}^{\ep}_{1,4}, \mathcal{M}^{\ep}_{1,6}$ and performing change of variables gives the same bounds.

To control $\mathcal{M}^{\ep}_{1,2}$ and $\mathcal{M}^{\ep}_{1,5}$ amounts to controlling $g^{\ep}$, which in turn requires the estimate in Proposition~\ref{prop:4.3} with $f^{\ep}$ replaced by $g^{\ep}$. We leave the proof to appendix~\ref{apx:pf_Mg_small}, but summarize the results in the following lemma.
\begin{lemma} \label{lem:M_geps_small} 
There exists some positive constant $c_8$ such that
\begin{equation}
|\mathcal{M}^{\ep}_{1,2} + \mathcal{M}^{\ep}_{1,5}| \leq c_8 \theta (\theta_{\mu} \theta_{\omega} \theta_{t})^{-\frac{1}{q}} \ep^{-(1+\frac{1}{q})}\,.
\end{equation}
\end{lemma}
This completes the proof for Theorem~\ref{them:measurem_decomp}.

%%%%%%%%%

\section{Numerical results}\label{sec:numerics}
Numerical evidence confirms our theoretical finding and verifies that stability deteriorates in the small Knudsen number regime.

\begin{figure}[htbp]
     \centering
     \begin{subfigure}[b]{0.4\textwidth}
         \centering
         \includegraphics[width=\textwidth]{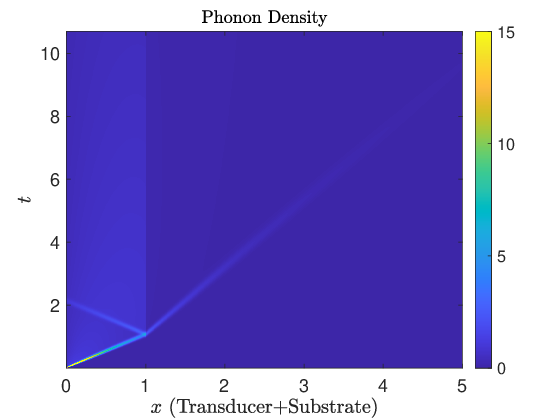}
         \caption{$\ep=1$.}
         \label{fig:density_fg_Kn1}
     \end{subfigure}
     \hfill
     \begin{subfigure}[b]{0.4\textwidth}
         \centering
         \includegraphics[width=\textwidth]{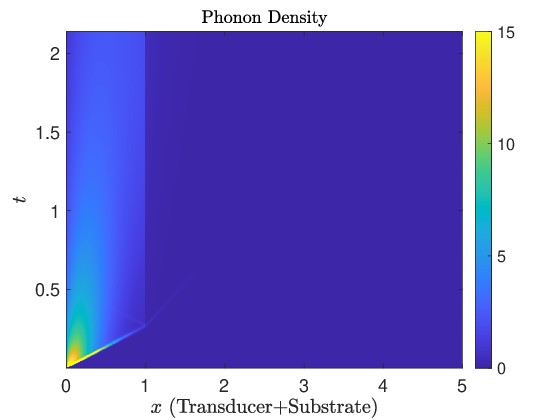}
         \caption{$\ep=0.25$.}
         \label{fig:density_fg_Kn025}
     \end{subfigure}
        \caption{Phonon density $\langle f^\ep\rangle(x) (0\leq x \leq 1)$, $\langle g^\ep\rangle(x) (1\leq x \leq 5)$ with different $\ep$ values. In this simulation of~\eqref{eq:forward_PTE_f_eps_couple} and~\eqref{eq:forward_PTE_g_eps_simple}, we use the following parameters: $\nu_\tr(\omega) = \omega, \tau_\tr(\omega)=1/\omega, \eta_{\tr} = (\tanh(10(\omega-1.5))-\tanh(2(\omega-1)))/4 + 0.9$.}
        \label{fig:density_fg_Kn}
\end{figure}

To conduct numerical experiments, we consider only system of $f^\ep$~\eqref{eq:forward_PTE_f_eps_couple}, disregarding $g^\ep$. This is a valid numerical approximation and is confirmed by Figure~\ref{fig:density_fg_Kn}. In this plot we show the density of $f^\ep, g^\ep$ (i.e., $\langle f^\ep\rangle$ and $\langle g^\ep\rangle$) with $f^\ep,g^\ep$ solving the coupled system~\eqref{eq:forward_PTE_f_eps_couple}-\eqref{eq:forward_PTE_g_eps_simple}. For both values of $\ep$, we observe that the phonon density of the substrate is much smaller than that of the transducer, and the transmitted phonon energy from transducer to substrate hardly bounce back to affect the transducer. In a sense $g^\ep$ appears as a downstream quantity affected by $f^\ep$, and that $f^\ep$ computation is essentially decoupled from the rest of the system.

To numerically verify the instability in small $\ep$ regime. We choose the following parameters: 
$\tau(\omega)=1/\omega, \nu(\omega)=\omega$ in~\eqref{eq:forward_PTE_f_eps_couple} in the examples below, and the source as
\begin{equation}\label{eq:source_phi_deltas}
\phi_i(x)=\frac{1}{\Delta \mu\Delta\omega\Delta t}\mathbf{1}_{\mu=\mu_0}\mathbf{1}_{\omega=\omega_i}\mathbf{1}_{t=0}\,,
\end{equation}
where $\mu_0=0.935$ and $\omega_i=\omega_{\min}+(i-1)\Delta\omega$ denotes the $i$-th grid of $\omega$. We use $\Delta\omega=0.05$ and $\omega_{\min}= \Delta\omega$ and $\omega_{max}=2$ such that there are 40 grid points. For $\mu$ discretization, we use $\Delta\mu = 0.01$ so that there are 200 grid points. In $x$ direction, we use $\Delta x = \min\{0.004,\ep/125\}$ and set $x\in[0,0.5]$ so that there are at least 126 grid points and for time we set $\Delta t=\min\{\ep \Delta x/2, \ep^2\}$ and use a variable stopping time $T=2\ep /2\mu_0\omega_i$. This is chosen to reflect the center of the measurement operator in~\eqref{eqn:measurement_delta}. The numerical scheme\footnote{The code is available~\url{https://github.com/Peiyi-wisc/PTE_stability}.} is upwinding in space, trapezoidal rule in $\mu, \omega$ and forward Euler in time, with $\Delta t$ chosen to satisfy CFL condition.

Moreover, we consider the following parametrization of $\eta(\omega)$:
\begin{equation}\label{eq:param_eta_omega}
\eta(\omega)=0.25\tanh{(10(\omega-a))} - 0.25\tanh{(2(\omega-b))}+0.5\,,
\end{equation}
so the stability will be shown only on these two parameters. The ground truth parameters are set to be $(a^\ast,b^\ast)=(1.5,1)$,

We first examine the dependence of $\ep$ of the cost function~\eqref{eq:L_eta_PDE_optim}. The data $\{d_{i,j}\}$ are generated by the ground truth $\eta^\ast(\omega)$.

\begin{figure}[htbp]
    \centering
    \includegraphics[scale = 0.3]{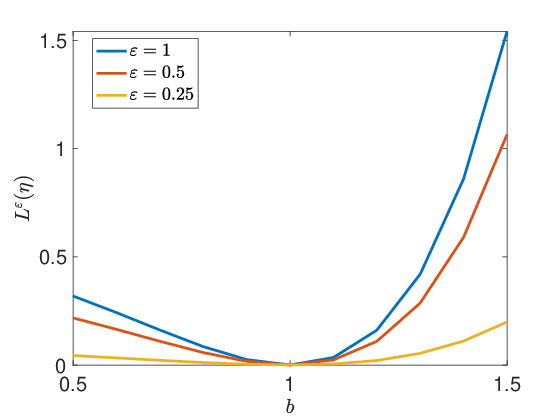}
    \caption{Cost function $L^\ep$ for different $\ep$ for varying $\eta(\omega)$.}
    \label{fig:cost_fn}
\end{figure}

In Figure~\ref{fig:cost_fn}, we plot the landscape of the loss function in $b\in[0.5, 1.5]$ for fixed $a=a^\ast$. It can be observed that as $\ep$ gets smaller, the cost function becomes flatter. This indicates that, changing from the ballistic regime to the diffusive regime, the temperature measurement error cannot distinguish different $\eta$, as shown in Corollary~\ref{thm:cor_loss_eta_lipschitz}.

\begin{figure}
    \centering
    \includegraphics[scale = 0.3]{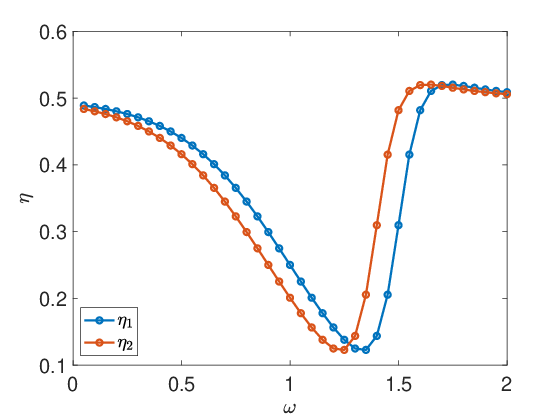}
    \caption{Two different forms of $\eta(\omega)$ used to compute the approximate values of $\mathcal{M}^\ep(\eta)$.}
    \label{fig:eta12}
\end{figure}

Next we examine the difference in $\mathcal{M}^\ep(\eta)$ when $\eta$ is perturbed. This is to run the PDE simulation for two different $\eta$'s (plotted in Figure~\ref{fig:eta12}), with $\eta_1 = \eta^\ast$, and $\eta_2$ is defined by~\eqref{eq:param_eta_omega} with $(a, b)=(1.4, 0.9)$.
% The $\eta_{1,2}$ are plotted in Figure~\ref{fig:eta12}, where 

\begin{figure}[htbp]
    \centering
    \includegraphics[scale = 0.2835]{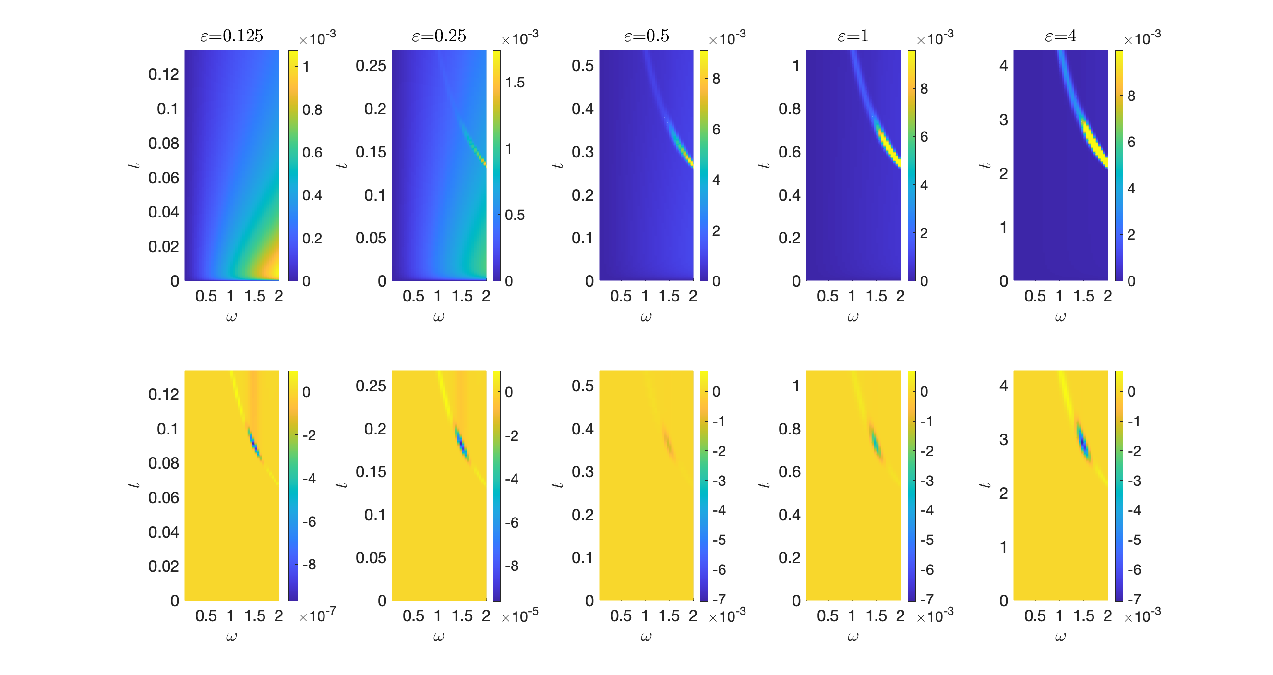}
    \caption{Profile of $\Lambda^\ep_{\eta_1}$ (first row) and the difference $\Lambda^\ep_{\eta_1}-\Lambda^\ep_{\eta_2}$ (second row) for varying $\ep$. The $x$-axis indicates the location of the source (in frequency) and the $y$-axis indicates the measurement in time. This plot shows that the sensitivity of the measurement map towards sources centered at various frequencies is severely reduced as $\ep$ decreases. Note the differences changes the order from $O(10^{-7}$ for small $\ep$ to $O(10^{-3})$ for big $\ep$.}
    \label{fig:A_eta12}
\end{figure}

In Figure~\ref{fig:A_eta12}, we plot both the measurement $\Lambda^\ep_{\eta_1}[\phi](t)$ (top row) and the difference $\Lambda^\ep_{\eta_1}[\phi](t) - \Lambda^\ep_{\eta_2}[\phi](t)$ (bottom row), for $\ep=0.125, 0.25, 0.5, 1, 4$ in different columns. Each plot is a function of $\omega$ and $t$. As $\ep$ increases, the differences $\Lambda^\ep_{\eta_1}[\phi](t) - \Lambda^\ep_{\eta_2}[\phi](t)$ becomes more enhanced, as shown by the order of the difference, changing from $O(10^{-7})$ to $O(10^{-5})$ to $O(10^{-3})$. The comparison between the last three columns is reflected by the colors of the difference.

\begin{figure}[htbp]
    \centering
    \begin{subfigure}[b]{0.32\textwidth} 
        \centering
        \includegraphics[width=\textwidth]{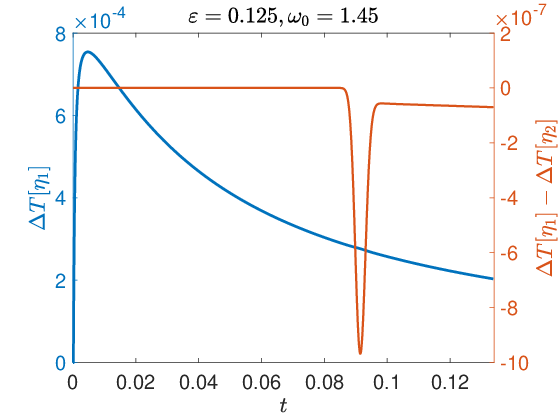}
        \caption{ }
        \label{fig:subfigA}
    \end{subfigure}
    \hfill 
    \begin{subfigure}[b]{0.32\textwidth}
        \centering
        \includegraphics[width=\textwidth]{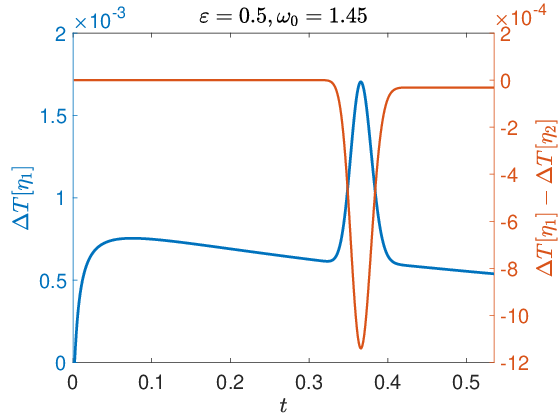}
        \caption{ }
        \label{fig:subfigB}
    \end{subfigure}
    \hfill 
    \begin{subfigure}[b]{0.32\textwidth}
        \centering
        \includegraphics[width=\textwidth]{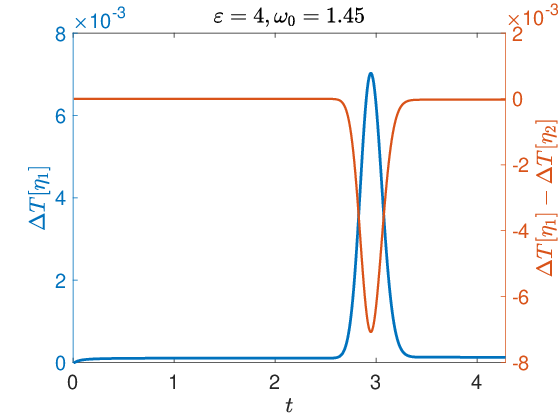}
        \caption{ }
        \label{fig:subfigC}
    \end{subfigure}
    \caption{Temperature reading in time, $\Delta T$, for $\eta_1$ and the temperature difference between $\eta_1$ and $\eta_2$, $\Delta T[\eta_1]-\Delta T[\eta_2]$, for $\ep$ being $0.125, 0.5, 4$.}
    \label{fig:temperature_time_reading}
\end{figure}

In Figure~\ref{fig:temperature_time_reading} we show one slide of Figure~\ref{fig:A_eta12} by slicing the $\omega=1.45$. The three plots respectively show the temperature reading, $\Delta T[\eta_1]$ in time, and the difference of $\Delta T[\eta_1]- \Delta T[\eta_2]$ for $\ep=0.125, \ep=0.5$ and $\ep=2$. It is clearly shown that when $\ep=4$, the temperature reading concentrates around $t=3$, which corresponds to the reflective behavior of phonons in this ballistic regime. However, when $\ep=0.125$, the temperature reading gradually decreases after the heat source is injected. The magnitude of the absolute difference $|\Delta T[\eta_1]- \Delta T[\eta_2]|$ in this case is about $10^{-7}$, and grows as the Knudsen number becomes larger, i.e., changing from the diffusive to ballistic regime. 

\begin{figure}[htbp]
    \centering
    \includegraphics[scale = 0.3]{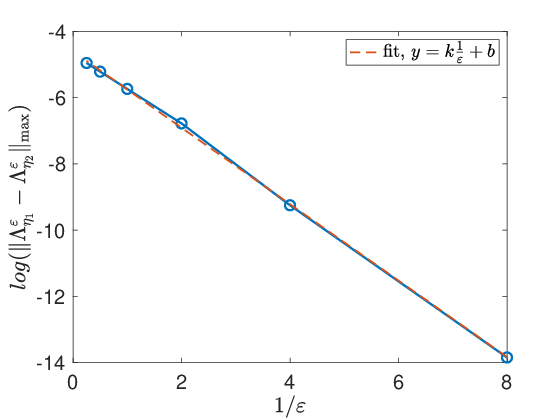}
    \caption{Maximum value of the difference between $\Lambda^{\ep}_{\eta_1}$ and $\Lambda^{\ep}_{\eta_2}$ in log scale for different $\ep$ values.}
    \label{fig:regression}
\end{figure}

In Figure~\ref{fig:regression} we present the dependence of $\|\Lambda^{\ep}_{\eta_1}-\Lambda^{\ep}_{\eta_2}\|_{\max}$ on $1/\ep$ in log scale. This shows the differences between two operators exponentially decrease with respect to $1/\ep$, confirming the stability estimate shown in Theorem~\ref{thm:stab_estim}.
% is also numerically justified. 
% For the given two reflective coefficients $\eta_{1,2}$, we compute the temperature $\Delta T(t,x=0)$ with source $\{\phi_i\}_{i=1}^{40}$~\eqref{eq:source_phi_deltas} respectively. 
% Therefore, at the discrete level, we compute the difference between the discretization of the two source-to-measurement operators~\eqref{eq:forward_operator_M_eta}, $\Lambda^{\eps}_{\eta_1}, \Lambda^{\eps}_{\eta_2}\in\mathbb{R}^{N_t\times 40} (N_t := \lceil \frac{T}{\Delta t}\rceil)$, and measure it in the maximum norm, i.e.,
% $\|\Lambda^{\eps}_{\eta_1}-\Lambda^{\eps}_{\eta_2}\|_{\max}\,.$
% This is plotted in the logarithmic scale against $\frac{1}{\eps}$ 

\section*{Use of Generative-AI Tools Declaration}
The authors declare they have used Artificial Intelligence (AI) tools in the creation of this article.
The authors used ChatGPT to assist with language editing and structure formulation.

%The acknowledgments section should not be numbered.
\section*{Acknowledgments}
P.~Chen and Q.~Li are supported by the National Science Foundation-DMS No.2308440. P.~Chen is also supported by the Hirschfelder Scholar Fellowship award with funding from the Wisconsin Foundation \& Alumni Association and the Department of Mathematics at the University of Wisconsin-Madison. 

This material is based upon work supported by the National Science Foundation under Grant No. DMS-2424139, while the authors were in residence at the Simons Laufer Mathematical Sciences Institute in Berkeley, California, during the Fall 2025 semester.

\appendix

\section{Model reduction}\label{apx:reduction}

In this section, we present the derivation of the coupling between $(\eta_{\tr},\zeta_{\tr})$ and $(\eta_{\sub}, \zeta_{\sub})$~\eqref{eq:tr_sub_couple}. This coupling relation reduces the four unknown coefficients at the interface to $(\eta_{\tr}, \zeta_{\tr})$. Together with the compatibility condition~\eqref{eq:eta_1_compatible_condition}, it is sufficient to determine only $\eta_{\tr}$ to fully solve the coupled system.

Across the interface $x=L_{\Tr}$, there should be no net flux for each frequency, i.e.,
\begin{equation}\label{eq:nonlinear_netflux}
\begin{aligned}
    & \int_0^1 \mu \nu_{\Tr}(\omega)F(x=L_{\Tr})\,\rd\mu + \int_{-1}^{0} (-\mu) \nu_{\Sub}(\omega)G(x=L_{\Tr})\,\rd\mu \\
    = & \int_{-1}^0 (-\mu) \nu_{\Tr}(\omega)F(x=L_{\Tr})\,\rd\mu + \int_{0}^{1} \mu \nu_{\Sub}(\omega)G(x=L_{\Tr})\,\rd\mu\,.
\end{aligned}
\end{equation}
The boundary conditions for $F,G$ at the interface $x=L_{\Tr}$ in~\eqref{eq:nonlinear-PTE-F} and~\eqref{eq:nonlinear-PTE-G} are then applied to compute the outgoing parts (RHS of~\eqref{eq:nonlinear_netflux}).
\begin{equation}
\begin{aligned}
    & \int_{-1}^0 (-\mu) \nu_{\Tr}(\omega) F(x=L_{\Tr})\,\rd\mu \\
    = & \, 
     \nu_{\Tr}F^\ast_0 (1-\eta_{\tr}-\zeta_{\tr}) \int_{-1}^0 (-\mu)\,\rd\mu
    + \nu_{\Tr} \eta_{\tr} \int_{-1}^0 (-\mu) F(x=L_{\Tr},-\mu)\,\rd\mu \\
    &\, + \nu_{\Tr} \zeta_{\tr} \int_{-1}^0 (-\mu)  G(x=L_{\Tr})\,\rd\mu\,,\\
    = & \, \frac{1}{2}\nu_{\Tr}F^\ast_0 (1-\eta_{\tr}-\zeta_{\tr}) 
    + \nu_{\Tr} \eta_{\tr} \int_{0}^1 \mu F(x=L_{\Tr})\,\rd\mu
    + \nu_{\Tr} \zeta_{\tr} \int_{-1}^0 (-\mu)  G(x=L_{\Tr})\,\rd\mu\,.
\end{aligned}
\end{equation}
Similarly, one has
\begin{equation}
\begin{aligned}
\int_{0}^{1} \mu \nu_{\Sub} G(x=L_{\Tr})\,\rd\mu 
= & \frac{1}{2} \nu_{\Sub} F^\ast_0(1-\eta_{\sub} - \zeta_{\sub}) \\
& + \nu_{\Sub} \eta_{\sub}\int_{-1}^{0} (-\mu) G(x=L_{\Tr}) \,\rd\mu 
+ \nu_{\Sub} \zeta_{\sub} \int_{0}^{1} \mu  F(x=L_{\Tr}) \,\rd\mu \,.
\end{aligned}
\end{equation}
We can group the same terms and obtain the following
\begin{equation}
\begin{aligned}
   0 = & (\nu_{\Tr} - \nu_{\Tr}\eta_{\tr} - \nu_{\Sub} \zeta_{\sub}) \int_0^1 \mu F(x=L_{\Tr})\,\rd\mu 
   + (\nu_{\Sub} - \nu_{\Sub} \eta_{\sub} - \nu_{\Tr} \zeta_{\tr} ) \int_{-1}^{0} (-\mu) G(x=L_{\Tr})\,\rd\mu \\
    & - \frac{1}{2}F^\ast_0 (\nu_{\Tr}(1-\eta_{\tr}-\zeta_{\tr}) + \nu_{\Sub}(1-\eta_{\sub} - \zeta_{\sub})) \,.
\end{aligned}
\end{equation}
Therefore, for this equality to hold any solution $(F,G)$, it recovers the condition~\eqref{eq:tr_sub_couple}. Together with the compatibility condition~\eqref{eq:eta_1_compatible_condition}, we can express the other coefficients at the interface in terms of the reflection coefficient $\eta_{\tr}$:  
\begin{equation}
\begin{aligned}
\zeta_{\tr} = \frac{1}{c}(1 - \eta_{\tr})\,,\ 
    \zeta_{\sub} = \frac{\nu_{\Tr}}{\nu_{\Sub}}(1-\eta_{\tr})\,, \
    \eta_{\sub} = 1 - \frac{\nu_{\Tr}}{\nu_{\Sub}} \frac{1}{c}(1 - \eta_{\tr})\,,
\end{aligned}
\end{equation}
and the same relations hold with nomarlized $\nu_{\tr}$ and $\nu_{\sub}$.

\section{Proof of Proposition~\ref{prop:charac}}\label{subsec:charac}
\begin{proof}
Using the method of characteristics, set
\begin{equation}
    \frac{\rd X}{\rd t}=\frac{\mu\nu}{\ep},\quad f_0^\ep=f_0^\ep(t)=f_0^\ep(t,X(t))\,.
\end{equation}
Then along a fixed characteristic with slope $\frac{\mu\nu}{\ep}$, 
\begin{equation}
    X(t)=X(t_0)+\frac{\mu\nu}{\ep}(t-t_0)\,,
\end{equation}
and we have
\begin{equation}\label{eqn:g_0_mu>0}
f_0^\ep(t,X(t))=f_0^\ep(t,X(t_0)+\frac{\mu \nu}{\ep}(t-t_0))e^{-\frac{t-t_0}{\ep^2\tau}}\,.
\end{equation}
For $\mu>0$, if the characteristic hits the boundary at $x=0$, we have $X(t_0)=0$ and
\begin{equation}
    t_0= t-\frac{X(t)}{\mu\nu/\ep}\,.
\end{equation}
This gives us
\begin{equation}\label{eqn:f_0^eps_mupos_phi}
f_0^\ep(t,X(t))=\phi\left(t-\frac{X(t)}{\mu\nu/\ep}, \mu, \omega \right)e^{-\frac{X(t)}{\mu\nu\tau\ep}}\,,
\quad t\geq \frac{X(t)}{\mu\nu/\ep}\,.
\end{equation}
If the trajectory hits initial data, we have $t-\frac{X(t)}{\mu\nu/\ep}<0$ and $f_0^\ep=0$.

Similarly, for $\mu<0$, if the trajectory hits boundary data, $X(t_0)=1$ and $t_0=t-\frac{1-X(t)}{-\mu\nu/ \ep}>0$. This gives $f_0^ \ep(t,X(t))=\eta(\omega)f_0^ \ep\left(t_0,1,\frac{1-X(t)}{|\mu|\nu(\omega)},-\mu,\omega\right)$ by the reflection boundary condition. It can be further traced back to the initial condition, using~\eqref{eqn:f_0^eps_mupos_phi}, that is,
\begin{equation}\label{eqn:g_0_mu<0}
f_0^ \ep(t,X(t))=\eta(\omega)\phi\left(t-\frac{2-X(t)}{-\mu\nu/ \ep}, -\mu, \omega \right)e^{\frac{2-X(t)}{-\mu\nu\tau \ep}}\,,
\end{equation}
which hits the boundary source $\phi$ at $x=0$ if $\frac{2-X(t)}{-\mu\nu\tau \ep}\geq 0$. If the trajectory hits initial data, we have $f_0^\ep(t_0)=0$.

Combining~\eqref{eqn:g_0_mu>0} and~\eqref{eqn:g_0_mu<0} gives
\begin{equation}
\begin{aligned}
     f_0^\ep(t,x,\mu,\omega) & =\begin{cases}
        \phi\big(t-\frac{\ep x}{\mu\nu},\mu,\omega\big)e^{-\frac{x}{\mu\nu\tau\ep}}, &\mu>0, t>\frac{\ep x}{\mu\nu}\\
        \eta(\omega)\phi\big(t-\frac{\ep(2-x)}{|\mu|\nu},|\mu|,\omega\big)e^{-\frac{(2-x)}{|\mu|\nu\tau\ep}},  &\mu<0, t>\frac{\ep(2-x)}{|\mu|\nu}>0\\
        0, & \text{otherwise}
    \end{cases}
\end{aligned}
\end{equation}

Similarly, for $f_1^\ep$, we have 
\begin{equation}
    f_1^\ep(t,X(t)) = f_1^\ep(t_0,X(t_0))e^{-\frac{t-t_0}{\ep^2\tau}} 
    + \int\limits_{t_0}^t\frac{1}{\ep^2\tau} \frac{C_\omega}{C_\tau}\langle f^\ep/\tau \rangle\left(s,\frac{\mu\nu}{\ep}(s-t_0)+ X(t_0)\right)e^{-\frac{t-s}{\ep^2\tau}}\,\rd s\,.
\end{equation}
For $\mu>0$, if the trajectory hits the boundary data at $x=0$, $t_0=t-\frac{X(t)}{\mu\nu/\ep}>0$, i.e., $t>\frac{X(t)}{\mu\nu/\ep}$, then we have
\begin{equation}\label{eq:f1_eps_boundary}
f_1^\ep(t,X(t)) = \int\limits_{t-\frac{X(t)}{\mu\nu/\ep}}^t \dfrac{1}{\ep^2\tau} \dfrac{C_\omega}{C_\tau}\langle f^\ep/\tau \rangle \left(s,\frac{\mu\nu}{\ep}(s-t)+X(t)\right)e^{-\frac{t-s}{\ep^2\tau}}\,\rd s\,.
\end{equation}
If the trajectory hits the initial data, we have $t_0=0$ and
\begin{equation}\label{eq:f1_eps_initial}
f_1^\ep(t,X(t)) = \int\limits_{0}^t \frac{1}{\ep^2\tau} \frac{C_\omega}{C_\tau} \langle f^\ep/\tau \rangle\left(s, \frac{\mu\nu}{\ep}(s-t)+ X(t)\right)e^{-\frac{t-s}{\ep^2\tau}}\,\rd s\,.
\end{equation}
Similarly for $\mu<0$, if the trajectory hits the boundary data at $x=1$, then $t_0=t-\frac{1-X(t)}{|\mu|\nu/\ep}>0$, i.e, $t>\frac{1-X(t)}{|\mu|\nu/\ep}$,
\begin{equation}
\begin{aligned}
f_1^\ep(t,X(t)) = & \,
\eta(\omega) f^\ep_{{1,\mu>0}}\left(t-\frac{1-X(t)}{|\mu|\nu/\ep},1\right) e^{-\frac{1-X(t)}{|\mu|\nu\ep\tau}} \\
& + \zeta(\omega) g^\ep_{{\mu<0}}\left(t-\frac{1-X(t)}{|\mu|\nu/\ep},1\right) e^{-\frac{1-X(t)}{|\mu|\nu\ep\tau}} \\
& + \int\limits_{t-\frac{1-X(t)}{|\mu|\nu/\ep}}^t \dfrac{1}{\ep^2\tau} \dfrac{C_\omega}{C_\tau} \langle f^\ep / \tau\rangle \left(s,\frac{\mu\nu}{\ep} (s-t) + X(t)\right)e^{-\frac{t-s}{\ep^2\tau}}\,\rd s\,.
\end{aligned}
\end{equation}
Using the formula for $f_1^\ep$ for $\mu>0$ gives that 
\begin{itemize}
    \item For $\dfrac{1-x}{\mu\nu/\ep}<t<\dfrac{2-x}{\mu\nu/\ep}$,
    \begin{equation}
    \begin{aligned}
        f^\ep_{1,\mu>0}\left(t-\frac{1-x}{\mu\nu/\ep}, 1\right) = & \int\limits_{0}^{t-\frac{1-x}{\mu\nu/\ep}}  \frac{1}{\ep^2\tau} \frac{C_\omega}{C_\tau} \langle f^\ep/\tau \rangle \left(s, 1+ \frac{\mu\nu}{\ep}\left(s-t+\frac{1-x}{\mu\nu/\ep}\right)\right)\\
        & \hspace{2in} \times e^{-\frac{1}{\ep^2\tau}(t-\frac{1-x}{\mu\nu/\ep}-s)}\,\rd s\,.
    \end{aligned}
    \end{equation}

    \item For $t>\dfrac{2-x}{|\mu|\nu/\ep}$,
    \begin{equation}
    \begin{aligned}
f^\ep_{1,\mu>0}\left(t-\frac{1-x}{\mu\nu/\ep}, 1\right) = &
\int\limits_{t-\frac{2-x}{\mu\nu/\ep}}^{t-\frac{1-x}{\mu\nu/\ep}} \frac{1}{\ep^2\tau} \frac{C_\omega}{C_\tau}\langle f^\ep /\tau \rangle \left(s, 1+\frac{\mu\nu}{\ep} \left(s-t + \frac{1-x}{\mu\nu/\ep}\right)\right) \\
& \hspace{2in} \times e^{-\frac{1}{\ep^2\tau}(t-\frac{1-x}{|\mu|\nu/\ep}-s)}\, \rd s\,,
\end{aligned}
\end{equation}
\end{itemize}
where the two cases of $X(t)$ hitting the boundary data~\eqref{eq:f1_eps_boundary} and the initial data~\eqref{eq:f1_eps_initial} are distinguished.
Then we apply the change of variable that $s':=t-s-\frac{1-x}{\mu\nu/\ep}$, and obtain that
\begin{itemize}
    \item For $\frac{1-x}{\mu\nu/\ep}<t<\frac{2-x}{\mu\nu/\ep}$
    \begin{equation}
        f^\ep_{1,\mu>0}\left(t-\frac{1-x}{\mu\nu/\ep},1\right)=
     \int\limits_{0}^{t-\frac{1-x}{\mu\nu/\ep}} \frac{1}{\ep^2\tau} \frac{C_\omega}{C_\tau} \langle f^\ep/\tau \rangle \left(t-s'-\frac{1-x}{\mu\nu/\ep}, 1-\frac{\mu\nu}{\ep}s'\right) e^{-\frac{s'}{\ep^2\tau}}\,\rd s'\,,
    \end{equation}

    \item For $t>\frac{2-x}{|\mu|\nu/\ep}$,
    \begin{equation}
    f^\ep_{1,\mu>0}\left(t-\frac{1-x}{\mu\nu/\ep},1\right)=\int\limits_{0}^{\frac{1}{\mu \nu/\ep}}\frac{1}{\ep^2\tau} \frac{C_\omega}{C_\tau} \langle f^\ep/\tau \rangle \left(t-s'-\frac{1-x}{|\mu|\nu/\ep}, 1 - \frac{\mu\nu}{\ep}s'\right) e^{-\frac{s'}{\ep^2\tau}}\,\rd s'\,.
\end{equation}
    
\end{itemize}

If the trajectory $X(t)$ hits the initial data, $t_0=0$, that is, $t<\frac{1-X(t)}{|\mu|\nu/\ep}$, then
\begin{equation}
    f_{1,\mu<0}^\ep(t,X(t))=\int\limits_{0}^{t} \frac{1}{\ep^2\tau} \frac{C_\omega}{C_\tau}\langle f^\ep/\tau \rangle\left(t-s, X(t)-\frac{\mu\nu}{\ep}s\right)e^{-\frac{s}{\ep^2\tau}}\,\rd s\,.
\end{equation}
Finally, we have that for $\mu<0$, the solution reads as the followings.
\begin{itemize}
    \item For $0<t\le\frac{\ep(1-x)}{|\mu|\nu}$,
    \begin{equation}
        f_1^\ep(t,x,\mu<0,\omega)=
        \frac{1}{\ep^2\tau} \frac{C_\omega}{C_\tau} \int\limits_0^{t}\langle f^\ep/\tau \rangle\left(t-s,x-\frac{\mu\nu }{\ep} s\right)e^{-\frac{s}{\ep^2\tau}}\,\rd s\,.
    \end{equation}

    \item For $\frac{\ep(1-x)}{|\mu|\nu} < t \leq \frac{\ep(2-x)}{|\mu|\nu}$,
    \begin{equation}
        \begin{aligned}
        & f_1^\ep(t,x,\mu<0,\omega) \\
        = \,
        & \frac{\eta(\omega)}{\ep^2\tau} \frac{C_\omega}{C_\tau} e^{-\frac{(1-x)}{\ep|\mu|\nu\tau}} \int\limits_0^{t-\frac{\ep(1-x)}{|\mu|\nu}}\langle f^\ep/\tau \rangle \left(t-s-\frac{\ep(1-x)}{|\mu|\nu}, 1+\frac{\mu\nu}{\ep}s\right) e^{-\frac{s}{\ep^2\tau}}\,\rd s\\
        & + \zeta(\omega) g^\ep_{{\mu<0}}\left(t-\frac{1-X(t)}{|\mu|\nu/\ep},1\right) e^{-\frac{1-X(t)}{|\mu|\nu\ep\tau}} \\
        & + \frac{1}{\ep^2\tau} \frac{C_\omega}{C_\tau} \int\limits_{0}^{\frac{\ep(1-x)}{|\mu|\nu}}\langle f^\ep/\tau \rangle \left(t-s, x -\frac{\mu\nu}{\ep}s \right)e^{-\frac{s}{\ep^2\tau}}\,\rd s\,.
        \end{aligned}
    \end{equation}

    \item For $t > \frac{\ep(2-x)}{|\mu|\nu}$,
    \begin{equation}
        \begin{aligned}
        & f_1^\ep(t,x,\mu<0,\omega) \\
        =\, & \eta(\omega)\frac{1}{\ep^2\tau} \frac{C_\omega}{C_\tau} e^{-\frac{(1-x)}{\ep|\mu|\nu\tau}} \int\limits_0^{\frac{\ep}{|\mu|\nu}}\langle f^\ep/\tau \rangle \left(t-s-\frac{\ep(1-x)}{|\mu|\nu}, 1+\frac{\mu\nu}{\ep}s\right) e^{-\frac{s}{\ep^2\tau}}\,\rd s\\
        & + \zeta(\omega) g^\ep_{{\mu<0}}\left(t-\frac{1-X(t)}{|\mu|\nu/\ep},1\right) e^{-\frac{1-X(t)}{|\mu|\nu\ep\tau}} \\
        & + \frac{1}{\ep^2\tau} \frac{C_\omega}{C_\tau} \int\limits_{0}^{\frac{\ep(1-x)}{|\mu|\nu}}\langle f^\ep/\tau \rangle \left(t-s, x-\frac{\mu\nu s}{\ep} \right)e^{-\frac{s}{\ep^2\tau}}\,\rd s\,.
        \end{aligned}
    \end{equation}
\end{itemize}
 
\end{proof}

\section{Proof of Proposition~\ref{prop:4.3}}\label{apx:scatter_theta_bound}
\begin{proof} 
We first change the variable with $\mu':=-\mu$, then let $y:=\frac{\mu'\nu s}{\ep}$, so that the integral is rewritten as
\begin{equation}
\begin{aligned}
\int\limits_{t_1-\theta}^{t_1+\theta} 
\int\limits_0^{\frac{\nu s}{\ep}} \int\limits_{\omega} \int\limits_{0}^t \psi_t\left(\frac{t-t_1}{\theta}\right)\frac{C_\omega}{\tau^2(\omega)} \frac{\ep}{\nu(\omega) s} \langle |f^\ep|/\tau \rangle\left(t-s, y\right) \mathbf{1}_{y<1}  \,\rd s\rd \omega \rd y \rd t\,,
\end{aligned}
\end{equation}
where we incorporate the fact that $\psi(t)$ has the compact support on $(-1,1)$.

By assumption~\ref{hyp1} that $\tau(\omega)$ is lower bounded: $\tau\ge \tau_{\min}$ and  H\"{o}lder inequality, 
\begin{equation}
    \begin{aligned}
      & \int\limits_{0}^{\frac{\nu s}{\ep}} \langle |f^\ep|/\tau \rangle(t-s,y)\mathbf{1}_{y<1}\,\rd y \\
      \le & \tau_{\min}^{\frac{1-q}{q}} \left(\int\limits_{0}^{\frac{\nu s}{\ep}}\int |f^\ep|^p(t-s,\cdot) C_\omega^{1-p} \mathbf{1}_{y<1}\, \rd\omega\rd\mu\rd y \right)^{\frac{1}{p}} 
      \left(\int\limits_{0}^{\frac{\nu s}{\ep}} \int C_\omega/\tau(\omega) \, \rd\omega\rd\mu\rd y\right)^{1/q}\,, \\
      &\le c_9 
      \|f^\ep\|_{L^p(C_\omega^{1-p}\,\rd y \rd\omega \rd\mu)} \left(\frac{\nu(\omega) s}{\ep}\right)^{\frac{1}{q}}\,,
    \end{aligned}
\end{equation}
where the constant $c_9(p,\tau_{\min},C_\tau) := \tau_{\min}^{-\frac{1}{p}} C_\tau^{1-\frac{1}{p}}$.
The above is then controlled by~\eqref{eq:f_g_Lp_estimate}, $\|f^\ep\|_{L^p(C_\omega^{1-p}\,\rd\omega \rd\mu\rd y)}\leq c_6 \|\phi\|_{L^p(0,t;L^p(C_\omega^{1-p}\, \rd\omega \rd\mu))}$. Moreover, by the definition of $\phi$~\eqref{eq:source_phi_t_mu_omega},
\begin{equation}\label{eq:phi_Lp_estimate}
\begin{aligned}
& \|\phi\|_{L^p(0,t;L^p(C_\omega^{1-p}\, \rd\omega \rd\mu))} 
\leq c_{10} (\theta_\mu\theta_\omega\theta_t)^{-1/q} \,,\\
& \hspace{1.2in} \text{with } c_{10} = \left(\int\limits_0^{\infty} \int\limits_{-1}^{1}
\int\limits_{0}^{\infty} \phi_t^p(s)\phi_{\mu}^p (\mu)\phi_{\omega}^p(\omega) C_{\omega}^{1-p}\,\rd s\rd\mu\rd\omega\right)^{1/p}\,.
\end{aligned}
\end{equation}
Applying this to~\eqref{eq:estimate_temp_feps} gives
\begin{equation}
    \begin{aligned}
    & \frac{1}{C_\tau^2} \frac{1}{\ep^2} 
        \int\limits_{0}^{t_{\max}} \int\limits_{-1}^0 \int\limits_{\omega} \int\limits_{0}^t\psi_t\left(\frac{t-t_1}{\theta}\right)\frac{C_\omega}{\tau^2(\omega)}\langle |f^\ep|/\tau \rangle\left(t-s,-\frac{\mu\nu s}{\ep}\right) \mathbf{1}_{\frac{|\mu|\nu s}{\ep}<1} \,\rd s\rd \omega \rd \mu \rd t \\
     \le &\, c_7 \ep^{-\left(1+\frac{1}{q}\right)} \theta (\theta_\mu\theta_\omega\theta_t)^{-\frac{1}{q}}  \,, 
    \end{aligned}
\end{equation}
where the constant $c_7=c_7(\phi,\psi,p,\tau_{\min},C_\tau)$ and can be computed by collecting the constants from the previous steps and Assumption~\ref{hyp1}.       
\end{proof}

\section{Proof of Lemma~\ref{lem:M_geps_small}}\label{apx:pf_Mg_small}
\begin{proof}
The singular decomposition technique can also be applied here. The solution $g^{\ep}$ to~\eqref{eq:forward_PTE_g_eps_couple} admits the decomposition that $g^{\ep}=:g_0^{\ep} + g_1^{\ep}$ such that 
\begin{equation}
\left\{
\begin{aligned}
\ep \partial_t g_0^{\ep} + \mu \nu(\omega)\partial_x g_0^{\ep} &= -\frac{1}{\ep}\frac{1}{\tau(\omega)} g_0^\ep\,, \hspace{2.05in} x\in[1, 1+L_{\sub}]\\
g_0^\ep (t,x=1,\mu,\omega) &= \eta(\omega)g_0^{\ep}(t,x=1,-\mu,\omega) + \zeta(\omega)f^{\ep}(t,x=1,\mu,\omega)\,,\hspace{.36in} \mu>0\\
g_0^\ep (t,x=L_s,\mu,\omega) &=0\,, \hspace{3.1in} \mu<0 \\
g_0^\ep (t=0,x,\mu,\omega) &=0\,, \hspace{3.1in} \mu<0
\end{aligned}
\right.
\end{equation}
here $\eta=\eta_{\sub},\zeta=\zeta_{\sub}$, and 
\begin{equation}
\left\{
\begin{aligned}
\ep \partial_t g_1^{\ep} + \mu \nu(\omega)\partial_x g_1^{\ep} &=\frac{1}{\ep}\frac{1}{\tau(\omega)} \left(\mathcal{L}g^\ep - g_1^\ep\right)\,, & x\in[1, 1+L_{\sub}]\\
g_1^\ep (t,x=1,\mu,\omega) &= 0\,, & \mu>0\\
g_1^\ep (t,x=L_s,\mu,\omega) &=0\,, & \mu<0 \\
g_1^\ep (t=0,x,\mu,\omega) &=0\,, & \mu<0
\end{aligned}
\right.
\end{equation}
where the subscripts of substrate on the parameters $\nu,\tau,\eta,$ and $\zeta$ are suppressed.
By the method of characteristics, 
\begin{equation*}
g_0^{\ep}(t,x=1,\mu<0,\omega) = 0\,, % 0\leq t \leq \frac{\ep (2 L_s)}{|\mu| \nu(\omega)}\,,
\end{equation*}
so $\mathcal{M}^{\ep}_{\mu<0}(g_0^{\ep})=0$, and
\begin{equation}\label{eq:M_f12_f15_g1}
\mathcal{M}^{\ep}_{1,2} + \mathcal{M}^{\ep}_{1,5} = 
\frac{2}{C_{\tau}} \int\limits_0^{t_{\max}} \psi(t)\int\limits_{-1}^0 \int \limits_{0}^{\infty} \zeta(\omega)
e^{-\frac{1}{\ep^2\tau} \frac{\ep}{|\mu|\nu}} g_1^{\ep}(t-\frac{\ep}{|\mu|\nu}, 1, \mu, \omega) \frac{1}{\tau(\omega)} \,\rd\omega\rd\mu \mathbf{1}_{t\geq \frac{\ep}{|\mu|\nu}} \rd t\,.
\end{equation}

Similarly, by the method of characteristics, we have:
\begin{itemize}
\item For $0\leq t\leq \frac{\ep L_s}{|\mu|\nu(\omega)}$, the solution hits the initial data, then
\begin{equation}
\begin{aligned}
g_1^{\ep}(t,1, \mu<0,\omega) 
& =
\int\limits_{0}^{t} \frac{1}{\ep^2\tau} e^{-\frac{t-s}{\ep^2 \tau}} \mathcal{L}[g^{\ep}](s, 1-\frac{\mu \nu(\omega)}{\ep} (t-s) ,\mu, \omega) \,\rd s\,,\\
& = \int\limits_{0}^{t} \frac{1}{\ep^2\tau} e^{-\frac{y}{\ep^2 \tau}} \mathcal{L}[g^{\ep}](t-y, 1-\frac{\mu \nu(\omega)}{\ep}y ,\mu, \omega) \,\rd y\,.
\end{aligned}
\end{equation}

\item For $t> \frac{\ep L_s}{|\mu|\nu(\omega)}$, the solution hits the boundary data, then
\begin{equation}
\begin{aligned}
g_1^{\ep}(t,1, \mu<0,\omega) 
= &
\int\limits_{t+\frac{\ep L_{\sub}}{\mu \nu(\omega)}}^{t} \frac{1}{\ep^2\tau} e^{-\frac{t-s}{\ep^2 \tau}}  \mathcal{L}[g^{\ep}](s, 1 - \frac{\mu\nu(\omega)}{\ep}(t-s),\mu, \omega) \,\rd s\,,\\
= & 
\int\limits_{0}^{\frac{\ep L_{\sub}}{|\mu| \nu(\omega)}} \frac{1}{\ep^2\tau} e^{-\frac{y}{\ep^2 \tau}}  \mathcal{L}[g^{\ep}](t-y, 1 - \frac{\mu\nu(\omega)}{\ep}y,\mu, \omega) \,\rd y\,,
\end{aligned}
\end{equation}
\end{itemize}
Note that for both cases, the following holds:
\begin{equation}
g_1^{\ep}(t,1, \mu<0,\omega) 
\leq 
\int\limits_{0}^{t} \frac{1}{\ep^2\tau} e^{-\frac{y}{\ep^2 \tau}}  \mathcal{L}[g^{\ep}](t-y, 1 - \frac{\mu\nu(\omega)}{\ep}y,\mu, \omega) \,\rd y\,.
\end{equation}
Plugging in the above to~\eqref{eq:M_f12_f15_g1} together with $\zeta(\omega)\leq 1$ gives:
\begin{equation}
\begin{aligned}
& |\mathcal{M}^{\ep}_{\mu<0}(f_{1,2}^{\ep} + f_{1,5}^{\ep})| \\
\leq & \,
\frac{2}{C_{\tau}} \int\limits_0^{t_{\max}} \psi(t)\int\limits_{-1}^0 \int \limits_{0}^{\infty} \int\limits_{0}^{t-\frac{\ep}{|\mu|\nu}} \frac{1}{\ep^2\tau^2} e^{-\frac{y}{\ep^2 \tau}}  \mathcal{L}[g^{\ep}](t-\frac{\ep}{|\mu|\nu}-y, 1 - \frac{\mu\nu}{\ep}y,\cdot) \,\rd y
\,\rd\omega\rd\mu \mathbf{1}_{t\geq \frac{\ep}{|\mu|\nu}} \rd t\,,
\\
\leq & \, 
\frac{2}{C_{\tau}} \int\limits_0^{t_{\max}} \psi(t)\int\limits_{-1}^0 \int\limits_{0}^{\infty} 
\frac{1}{\ep^2\tau^2} \int\limits_{0}^{t-\frac{\ep}{|\mu|\nu}} 
\frac{C_{\omega}}{C_{\tau}} \langle |g^{\ep}|/\tau\rangle(t-\frac{\ep}{|\mu|\nu}-y, 1 - \frac{\mu\nu}{\ep}y)\,\rd y \rd\omega\rd\mu \mathbf{1}_{t\geq \frac{\ep}{|\mu|\nu}} \rd t\,,
\end{aligned}
\end{equation}
Let $z:=y+\frac{\ep}{|\mu| \nu(\omega)}$, by this change of variable, we obtain
\begin{equation}
\begin{aligned}
\int\limits_{0}^{t-\frac{\ep}{|\mu|\nu}} 
\langle |g^{\ep}|/\tau\rangle(t-\frac{\ep}{|\mu|\nu}-y, 1 - \frac{\mu\nu(\omega)}{\ep}y)\, \rd y 
= 
\int\limits_{\frac{\ep}{|\mu|\nu}}^{t} 
\langle |g^{\ep}|/\tau\rangle(t-z, \frac{|\mu|\nu(\omega)}{\ep}z )\, \rd z\,,
\end{aligned}
\end{equation}
then let $\tilde{\mu}:=\frac{|\mu|\nu(\omega) z}{\ep}$, the measurement integral can be further rewritten as
\begin{equation}
\begin{aligned}
|\mathcal{M}^{\ep}_{1,2} + \mathcal{M}^{\ep}_{1,5}| 
& \leq \frac{1}{\ep}\frac{2}{C_{\tau}^2} \int\limits_0^{t_{\max}} \int\limits_{0}^{\infty}  
\frac{C_{\omega}}{\tau^2 \nu(\omega)} \, \rd\omega 
\int\limits_0^{\frac{\nu z}{\ep}}
\int\limits_{0}^{t} \langle |g^{\ep}|/\tau\rangle(t-z, \tilde{\mu}) 
\frac{1}{z}
\,\rd z\rd\tilde{\mu}  \mathbf{1}_{t\geq \frac{\ep}{|\mu|\nu}} \psi\,\rd t\,.
\end{aligned}
\end{equation}
For the above integral we can apply H\"older inequality and the bound in $L^p$-norm (lemma~\ref{lem:prop_31_Lp_norm} and the estimate~\eqref{eq:phi_Lp_estimate}) to obtain
\begin{equation}
\begin{aligned}
& |\mathcal{M}^{\ep}_{1,2} + \mathcal{M}^{\ep}_{1,5}| \\
\leq & \,\frac{1}{\ep}\frac{2}{C_{\tau}^2} \int\limits_0^{t_{\max}} \psi(t)\int\limits_{0}^{\infty}  
\frac{C_{\omega}}{\tau^2 \nu(\omega)} \, \rd\omega 
\int\limits_{0}^{t} 
\int\limits_0^{\frac{\nu(\omega)z}{\ep}}
\langle |g^{\ep}|/\tau\rangle(t-z, \tilde{\mu}) 
\frac{1}{z}
\,\rd\tilde{\mu} \rd z \mathbf{1}_{t\geq \frac{\ep}{|\mu|\nu}} \rd t\,,\\
\leq & \,\tau_{\min}^{\frac{1-q}{q}} \frac{1}{\ep^{1+1/q}}\frac{2 C_{\tau}^{1/q}}{C_{\tau}^2} \int\limits_0^{t_{\max}} \psi(t)\int\limits_{0}^{\infty}  
\frac{C_{\omega}}{\tau^2 \nu^{1/p}} \, \rd\omega 
\int\limits_{0}^{t} \frac{1}{z^{1/p}} \|g^{\ep}(t-z,\cdot)\|_{L^p(C_{\omega}^{1-p}\,\rd y\rd\omega\rd\mu)}
\mathbf{1}_{t\geq \frac{\ep}{|\mu|\nu}} \,\rd z\rd t\,,\\
\leq & \, c_8 \ep^{-(1+1/q)} \theta (\theta_\mu \theta_\omega \theta_t)^{-\frac{1}{q}}\,,
\end{aligned}
\end{equation}
where we recall assumption~\eqref{eq:assump_bdd_moment}.
\end{proof}

%%%%%%%%%%%%%%%%%%%%%%%%%%%%%%%%%%%%%%%%%%%%%%%%%%%%%%
%          7. REFERENCES SECTION
%%%%%%%%%%%%%%%%%%%%%%%%%%%%%%%%%%%%%%%%%%%%%%%%%%%%%%

%       READ THIS SECTION CAREFULLY

% Each of the references below MUST be cited in your article above. Do not include references that are not cited in your article.

% Follow the examples below carefully. We strongly suggest that you copy and paste your reference information directly into our examples.

% List all references in alphabetical order according to the first author's last name.

% Verify each URL works correctly and can be accessed properly. Your URL links should be to reputable websites. The command line for a website link begins with: \url{ }

% Do not add MR or DOI numbers to your references. AIMS production staff will add this information.

% Using BibTex is not recommended but can be handled.

\medskip
% The information below will be filled in by AIMS production staff.
Received xxxx 20xx; revised xxxx 20xx; early access xxxx 20xx.
\medskip

\end{document}